\documentclass[11pt,reqno,a4paper]{amsart}

\usepackage[margin=25mm]{geometry}

\usepackage{graphicx}
\usepackage{multirow}
\usepackage{amsmath,amssymb,amsfonts}
\usepackage{amsthm}
\usepackage{mathrsfs}
\usepackage[title]{appendix}
\usepackage{xcolor}
\usepackage{textcomp}
\usepackage{manyfoot}
\usepackage{booktabs}
\usepackage{listings}
\usepackage{mathtools}
\usepackage{placeins}

\newtheorem{theorem}{Theorem}
 
\newtheorem{lemma}{Lemma}

\theoremstyle{remark}
\newtheorem{example}{Example}
\newtheorem{remark}{Remark}

\theoremstyle{definition}
\newtheorem{definition}{Definition}

\DeclareMathOperator{\M}{\mathbb E}
\newcommand{\R}{\mathbb{R}}

\DeclareMathOperator{\cov}{\mathbf{cov}}
\DeclareMathOperator{\var}{\mathbf{Var}}
\allowdisplaybreaks

\begin{document}

\title[Parameter estimation in models with explosive drift and Gaussian noise]{Parameter estimation in models with explosive drift and~Gaussian noise}

\author{Yuliya Mishura}
\address[Y. Mishura]{Department of Probability Theory, Statistics and Actuarial Mathematics, Taras Shevchenko National University of Kyiv, Kyiv, Ukraine}
\email{yuliyamishura@knu.ua}

\author{Kostiantyn Ralchenko}
\address[K. Ralchenko]{Department of Probability Theory, Statistics and Actuarial Mathematics, Taras Shevchenko National University of Kyiv, Kyiv, Ukraine}
\address[K. Ralchenko]{School of Technology and Innovations, University of Vaasa, Vaasa, Finland}
\email{kostiantynralchenko@knu.ua}
\thanks{Kostiantyn Ralchenko acknowledges support from the
Research Council of Finland (decision no.~367468).}

\author{Mykyta Yakovliev}
\address[M. Yakovliev]{Department of Mathematical Analysis, Taras Shevchenko National University of Kyiv, Kyiv, Ukraine}
\email{mykyta.yakovliev@knu.ua}

\begin{abstract}
We study drift parameter estimation for a positive stochastic process with singular mean-reverting drift driven by Gaussian noise.
For Riemann--Liouville fractional Brownian noise with Hurst parameter $H>1/2$, an observable fractional transformation reduces the model to a semimartingale with drift linear in the unknown parameters.
Using the quasi-likelihood approach, we derive explicit estimators of both drift coefficients.
For a broader class of Gaussian driving processes, we construct a strongly consistent estimator of the ratio $a/b$, which has a natural mean-reversion interpretation, and verify the required assumptions for several Gaussian models.
The theoretical results are complemented by a simulation study investigating the effects of the observation horizon and the discretization step on the finite-sample accuracy of the proposed estimators.
\end{abstract}



\maketitle

\section{Introduction}\label{sec1}

We consider statistical inference for the positive solution of the
stochastic equation
\begin{equation}\label{eq:model}
X_t = X_0 + a\int_0^t\frac{ds}{X_s}
- b\int_0^tX_s\,ds + G_t,
\qquad t\ge0,
\end{equation}
where $X_0>0$, $a>0$, $b>0$, and $G=\{G_t,t\ge0\}$ is a Gaussian
driving process. Our main objective is the estimation of the drift
parameters $a$ and $b$ from continuous observation of $X$.

Equation \eqref{eq:model} is closely related to the
Cox--Ingersoll--Ross (CIR) model introduced in
\cite{CoxIngersollRoss1985}. In the fractional setting, its solution can
be interpreted as a square-root process associated with a
Gaussian-driven CIR-type model; see, in particular,
\cite{MishYur2,MishYur3,MishYur1}. Such extensions are motivated by
the use of Gaussian models with memory in continuous-time finance; see, for example, \cite{ComteRenault1996} and the references therein.
The singular term $a/X_t$ prevents the process from reaching zero under suitable regularity assumptions on the noise, whereas the term
$-bX_t$ produces mean reversion. In the pathwise setting of this paper, the ratio
$\theta=\frac{a}{b}$
has a natural interpretation as the equilibrium level of the associated
squared CIR-type process, while $\sqrt{a/b}$ is the equilibrium point
of the drift of $X$ itself.

From a computational point of view, equations of the form
\eqref{eq:model} are also of interest because the singular drift and
Gaussian noise with memory make their numerical treatment nontrivial.
Positivity-preserving approximation schemes for fractional CIR and, more generally, sandwiched stochastic equations have been developed in \cite{DiNunnoMishuraYT2023,HongHuangKamraniWang2020}; see also \cite{MishYur3} for a systematic treatment. They provide effective tools for reliable numerical simulation and for studying the pathwise and financial properties of these models.

For applications, however, numerical approximation of the dynamics is only one part of the problem: the drift coefficients entering the
equation are typically unknown and have to be inferred from observed
data. Thus numerical simulation and statistical calibration are
closely connected. In the present setting, this interaction is
particularly important, since the estimators themselves involve
numerical approximation of fractional operators and are therefore
affected by both the observation horizon and the discretization step.
Accordingly, alongside the theoretical estimation results, we study
numerically the finite-sample behavior of the proposed estimators and
the impact of discretization on their accuracy.

Stochastic equations with singular drift driven by fractional Brownian
motion have attracted considerable attention. In particular,
\cite{HuNualartSong2008} established existence and uniqueness results
for a class of singular equations driven by fractional Brownian motion
with Hurst parameter $H>1/2$ and discussed their connection with
fractional CIR models. For the particular square-root type equation
underlying \eqref{eq:model}, existence, uniqueness, and strict
positivity were studied in \cite{MishYur1}, and subsequently extended
to more general singular drifts and Gaussian driving processes in
\cite{MishYur2,MishYur3}.

Statistical inference for diffusion processes driven by standard
Brownian motion is by now a well-developed subject; see, for example,
\cite{Kutoyants2004}. For fractional and other non-semimartingale
models, however, parameter estimation is substantially more delicate.
A classical example is the fractional Ornstein--Uhlenbeck model:
maximum likelihood estimation based on an associated fundamental
martingale was developed in \cite{KleptsynaLeBreton2002},
least-squares estimation was studied in \cite{HuNualart2010}, and
minimum contrast estimation was considered in \cite{Bishwal2011}.
These results illustrate both the feasibility and the additional
difficulties of statistical inference under Gaussian noise with memory.
In comparison, considerably less is known for nonlinear equations with
singular drift such as \eqref{eq:model}.

A natural first approach to estimating $(a,b)$ is maximum likelihood.
For ordinary diffusions this typically relies on Girsanov's theorem.
For Gaussian models with memory, however, the corresponding
change-of-measure argument depends strongly on the structure of the
driving process and requires additional absolute-continuity and
exponential-integrability conditions. This becomes particularly
cumbersome in the present nonlinear two-parameter setting. We therefore
proceed differently.

Our first contribution concerns the case where $G$ is a
Riemann--Liouville fractional Brownian motion with known Hurst parameter
$H>1/2$. Its simple Volterra representation in terms of a standard
Wiener process allows us to apply an observable fractional
transformation to \eqref{eq:model}. The transformed process is a
Brownian semimartingale whose drift depends linearly on $a$ and $b$.
Using the martingale quasi-likelihood approach of
\cite{Heyde1994,Heyde1997}, we construct an explicit quasi-score and
derive closed-form estimators of both drift parameters. The
quasi-likelihood construction avoids the need to establish mutual
absolute continuity of the parameterized laws or additional
exponential-martingale conditions. At the same time, whenever the
ordinary Girsanov likelihood is well defined, the resulting estimator
coincides with the corresponding maximum likelihood estimator. We also
give equivalent scalar quasi-likelihood and minimum-contrast
interpretations of the estimator.

Our second contribution is of a different nature and applies to a much
broader class of Gaussian driving processes. Rather than estimating
$a$ and $b$ separately, we consider the parameter ratio
$\theta=a/b$, which has the mean-reversion interpretation described
above. Under local H\"older regularity of order greater than $1/2$ and
a mild sublinear-growth condition on the driving process, we prove that
\[
\widehat\theta_T
=
\frac{\int_0^T X_s\,ds}
     {\int_0^T X_s^{-1}\,ds}
\to
\frac{a}{b}
\qquad\text{a.s. as }T\to\infty.
\]
Thus $\widehat\theta_T$ is a strongly consistent estimator of the
mean-reversion parameter $\theta$. The result does not require a
particular Volterra representation of the noise. We verify its
assumptions for several important examples, including fractional
Brownian motion, Riemann--Liouville fractional Brownian motion,
Gaussian Volterra processes with power-type kernels, tempered
fractional Brownian motions, and multifractional Brownian motion.

The theoretical results are complemented by a simulation study
investigating the finite-sample performance of the proposed estimators
and their numerical implementation. In particular, we examine how the
discretization step and the observation horizon affect the bias and the
root mean square error of the quasi-likelihood estimators, as well as
the behavior of the ratio estimator under Riemann--Liouville and
standard fractional Brownian noise.

The remainder of the paper is organized as follows.
Section~\ref{sec:prelim} collects the necessary facts on Gaussian
processes, H\"older regularity, and existence and positivity of the
solution. Section~\ref{sec2} develops the quasi-likelihood estimation
procedure for Riemann--Liouville fractional Brownian noise.
Section~\ref{sec:ratio} establishes strong consistency of the estimator
of $a/b$ and verifies the assumptions for several classes of Gaussian
processes. Section~\ref{sec4} contains the simulation study and the
computational results.

\section{Preliminaries}
\label{sec:prelim}
We start with two closely related definitions, of fractional Brownian motion and Riemann--Liouville fractional Brownian motion. Note, however, that, as a rule, these definitions are given differently, namely, for fractional Brownian motion, its covariance function is given, and then its compact-interval representation is provided, while for the Riemann--Liouville fractional Brownian motion, oppositely, compact-interval representation comes the first. Partly, this is a historical tradition, partly, this is explained by the fact that compact-interval representation of fBm is more complicated while for the Riemann--Liouville fractional Brownian motion it is very simple. We preserve this tradition. 

Throughout the paper, all stochastic processes are defined on a
filtered probability space
$(\Omega,\mathcal F,\{\mathcal F_t\}_{t\ge0},\mathbb P)$
satisfying the usual conditions. The notation
$W=\{W_t,t\ge0\}$ stands for a standard
$\{\mathcal F_t\}$-Wiener process.

\begin{definition}
        A \emph{fractional Brownian motion} (fBm) with Hurst parameter $H\in(0,1)$ is a centered Gaussian process $W^H=\{W_t^H,\ t\ge 0\}$ with covariance function
    \begin{equation*}
        \mathbb{E}[W_t^H W_s^H] = \frac{1}{2} \left(  t^{2H} +  s^{2H} - |t-s|^{2H} \right), \qquad s,t\ge  0.
    \end{equation*}
    
    Besides this covariance characterization, according to \cite{NorrosValkeilaVirtamo1999}, fBm admits a compact-interval representation in terms of a standard  Wiener process and Volterra-type kernel  
    $W_t^H = \int_0^t z_H(t,s)dW_s$, $t\ge 0$,
    where, for $H>\frac12$ ,
    \begin{equation*}
        z_H(t,s) = c_H \left(H-\tfrac12\right) s^{\frac12-H} \int_s^t u^{H-\frac12}(u-s)^{H-\frac32}\,du, \qquad 0<s<t.
    \end{equation*}
    with
    $c_H = \left( \frac{ 2H\Gamma\left(\frac{3}{2}-H\right) }{ \Gamma\left(H+\frac{1}{2}\right)\Gamma(2-2H) } \right)^{1/2}$.
\end{definition}

\begin{definition}\label{def:rlfbm}
 A \emph{Riemann--Liouville fractional Brownian motion} with Hurst parameter $H\in(0,1)$ is the centered Gaussian process $B^H=\{B_t^H,t\geq 0\}$ defined by 
\begin{equation*}
    B_t^H = \int_0^t (t-s)^{H-\frac{1}{2}}\,dW_s, \qquad t\geq 0.
\end{equation*}

Its covariance function is given by 
\begin{equation*}
    \cov(B_t^H,B_s^H) = \mathbb{E}\bigl[B_t^H B_s^H\bigr] = \int_0^{t\wedge s} (t-u)^{H-\frac{1}{2}} (s-u)^{H-\frac{1}{2}}\,du, \qquad s,t\geq 0.
\end{equation*}
In particular, 
$\var(B_t^H)=\frac{t^{2H}}{2H}$, $t\geq 0$.
\end{definition} 

\begin{remark} Again, this is just a tradition not to normalize Riemann--Liouville fBm, and, since we assume that $H$ is known, it will be of no difference for us. 
\end{remark}
\begin{definition}
 A stochastic process $U=\{U_t,t\ge0\}$ is called \emph{$\lambda$-Hölder continuous} on the interval $[0,T]$, $\lambda\in(0,1]$, if there exists a random variable $\Lambda=\Lambda(T,\omega)>0$ such that
\begin{equation}\label{Holderss}
    |U_t-U_s| \le \Lambda |t-s|^\lambda, \qquad s,t\in[0,T].
\end{equation}
\end{definition}
\begin{remark} \label{rem:holder} It is very well known that both fBm $W^H$ and Riemann--Liouville fBm $B^H$ are H\"older continuous on any interval $[0,T]$ up to order $H$, i.e., are $\lambda$-H\"older continuous on any interval $[0,T]$ of any order $0<\lambda<H$. See, e.g., \cite[Examples 3.4 and 3.5]{MishYur3}.  Moreover, for both of these processes and for any fixed $T>0$ the random variable $\Lambda$ can be
chosen in the explicit form given in \cite[Example~3.2, eq.~(3.5)]{MishYur3}, and then $\mathbb{E}\Lambda^{r}<\infty$ for every $r>0$.
\end{remark}

Now consider the stochastic differential equation with explosive drift and Gaussian noise of the form 
\begin{equation}\label{eq:modela}
X_t = X_0
+ \int_0^t \frac{a}{(X_s)^{\frac{\alpha}{1-\alpha}}}\,ds
- b \int_0^t X_s\,ds
+ G_t,
\end{equation} where $\alpha\in(0,1)$ is some constant, $a>0$, $b\ge 0$, $X_0>0$, and $G_0=0$. 
Sufficient existence--uniqueness conditions for \eqref{eq:modela} are provided in \cite[Sect.~4.1, eq.~(4.1)]{MishYur3}. In the simplified form, this result is contained in the following theorem.

\begin{theorem}\label{theorem1} Let $G$ be a $\lambda$-H\"older continuous process on the interval $[0,T]$ for some $\lambda\in(0,1)$, and assume that $\alpha+\lambda>1$.
Then equation \eqref{eq:modela} has a unique strictly positive solution on $[0,T]$.
\end{theorem}

\begin{remark} For equation \eqref{eq:model}, which is the subject of our consideration, we have $\alpha=1/2$, so that $\frac{\alpha}{1-\alpha}=1$. Therefore, condition $\alpha+\lambda>1$ is reduced to $\lambda>\frac12$. If we consider fBm or Riemann--Liouville fBm with index $H>1/2$ as the noise, then, for any fixed $1/2<\lambda<H$, the sample paths of the driving process are $\lambda$-H\"older continuous. Therefore, in both cases equation \eqref{eq:model} has a unique strictly positive solution on any interval $[0,T]$. Note also that, contrary to our intuition, the bigger power $\frac{\alpha}{1-\alpha}$ is better. 
\end{remark}
\begin{remark}\label{remark4}
In order to establish  the H\"older continuity of  the solution $X$, suppose that $G$ is $\lambda$-H\"older continuous on
$[0,T]$ for some $\lambda>1/2$, and that the H\"older constant $\Lambda$
in \eqref{Holderss} satisfies
$\mathbb E\Lambda^r<\infty$, $r>0$.
For fBm and Riemann--Liouville fBm with Hurst parameter $H>1/2$,
these assumptions hold for every $\lambda\in(1/2,H)$ by
Remark~\ref{rem:holder}.

First, by Theorem~3.3 and formula~(3.39) in \cite{MishYur3}, there exist deterministic constants $C_1, C_2>0$, independent of the realization of the driving process, such that  a.s.,
\begin{equation}\label{infsup}
C_1(1+\Lambda)^{-\frac{1}{2\lambda-1}}
    \le \inf_{t\in[0,T]}X_t \le \sup_{t\in[0,T]}X_t \le
    C_2(1+\Lambda).
\end{equation}

Second, it then follows from formula~(3.40) and Theorem~3.4 in \cite{MishYur3} that
\begin{equation*}
    |X_t-X_s| \leq \Upsilon |t-s|^\lambda, \qquad s,t\in[0,T],
\end{equation*}
where $\Upsilon$ is a finite random H\"older constant of the solution. Its dependence on the random quantities can be bounded as follows: 
\begin{equation*}
    \Upsilon \leq   C_3 (1+\Lambda)^{1+\frac{2}{2\lambda-1}},
\end{equation*}
where $C_3$ is a deterministic constant. 
Therefore,
\begin{equation}\label{momentgamma}
\mathbb{E}\Upsilon^r<\infty, \quad r>0.  
\end{equation}

The relations \eqref{infsup} and \eqref{momentgamma} 
imply that 
\begin{equation}\label{eq:moments}
    \mathbb{E}\biggl[\sup_{t\in[0,T]}X_t^{\,r}\biggr]<\infty,\quad
    \text{and} \quad\mathbb{E}\biggl[\sup_{t\in[0,T]}X_t^{-r}\biggr]<\infty,
    \quad r > 0.
    \end{equation} 
\end{remark}

\section{Quasi-likelihood estimation of the drift parameters under
Riemann--Liouville fractional Brownian noise}\label{sec2}

In this section, we construct quasi-likelihood estimators for the drift
parameters of equation \eqref{eq:model} when the driving process is a
Riemann--Liouville fractional Brownian motion. Its representation in
terms of a standard Wiener process, given in
Definition~\ref{def:rlfbm}, allows us to apply an observable fractional
transformation to the original model and obtain a Brownian
semimartingale whose drift is linear in the unknown parameters. We then
apply the martingale quasi-likelihood approach of \cite{Heyde1994};
see also \cite[Sect.~8.3]{Heyde1997}, to derive an explicit estimating
function and the corresponding estimators.

An advantage of this approach is that it does not require establishing
mutual absolute continuity of the probability laws corresponding to
different parameter values or verifying additional
exponential-martingale conditions. As discussed in
Remark~\ref{rem:QLE-MLE}, whenever the ordinary likelihood is well
defined, the resulting quasi-likelihood estimator coincides with the
corresponding maximum likelihood estimator.

Let $B^H$ be a Riemann--Liouville fractional Brownian motion with
$H\in(\tfrac12,1)$, and consider
\begin{equation}\label{eq:model2}
X_t
=
X_0
+a\int_0^t\frac{ds}{X_s}
-b\int_0^tX_s\,ds
+\sigma B_t^H,
\qquad t\ge0,
\end{equation}
where $\sigma>0$ and $H\in(\tfrac12,1)$ are assumed to be known,
whereas the parameters $a>0$ and $b>0$ are unknown.

\subsection{Transformation to a Brownian-driven model}
\label{subsec:transformation}

We transform \eqref{eq:model2} by means of Rie\-mann--Liouville
fractional differentiation; see, for example, \cite{AMP2023}. For
$t>0$, define
\begin{gather}\label{eq:def-Y}
Y_t
\coloneqq
\frac{d}{dt}
\left(
\int_0^t(t-s)^{\frac12-H}X_s\,ds
\right)
=
\Gamma\left(\tfrac32-H\right)
\left(D_{0+}^{H-\frac12}X\right)(t),
\shortintertext{and}
\label{eq:def-Z}
Z_t
\coloneqq
\frac{d}{dt}
\left(
\int_0^t(t-s)^{\frac12-H}X_s^{-1}\,ds
\right)
=
\Gamma\left(\tfrac32-H\right)
\left(D_{0+}^{H-\frac12}X^{-1}\right)(t).
\end{gather}

Let us first verify that these fractional derivatives are well defined.
Fix $T>0$ and choose $\lambda\in(1/2,H)$. Since $B^H$ is
$\lambda$-H\"older continuous on $[0,T]$, the same is true for $X$.
Indeed, the strict positivity and continuity of $X$ imply that
\[
m_T
\coloneqq
\inf_{0\le t\le T}X_t>0
\qquad\text{a.s.},
\]
and hence the drift term in \eqref{eq:model2} is Lipschitz continuous
on $[0,T]$. Moreover,
\[
|X_t^{-1}-X_s^{-1}|
\le
m_T^{-2}|X_t-X_s|,
\]
so $X^{-1}$ is also $\lambda$-H\"older continuous.

Since $\lambda>H-\frac12$, the Riemann--Liouville fractional
derivatives of order $H-\frac12$ of both $X$ and $X^{-1}$ are well
defined and admit the Marchaud representations; see
\cite[Sections~13.1 and~13.4]{SMK}. Consequently,
\begin{gather}\label{forY}
Y_t
=
\frac{X_t}{t^{H-\frac12}}
+
\left(H-\tfrac12\right)
\int_0^t
\frac{X_t-X_s}{(t-s)^{H+\frac12}}\,ds,
\qquad t>0,
\shortintertext{and}
\label{forZ}
Z_t
=
\frac{X_t^{-1}}{t^{H-\frac12}}
+
\left(H-\tfrac12\right)
\int_0^t
\frac{X_t^{-1}-X_s^{-1}}
     {(t-s)^{H+\frac12}}\,ds,
\qquad t>0.
\end{gather}
Thus $Y$ and $Z$ are completely determined by the observed trajectory
of $X$ and are therefore observable transformations of the original
process.

We now rewrite \eqref{eq:model2} in terms of $Y$ and $Z$.

\begin{theorem}\label{thm:transformed-model}
Let $X$ be the solution of \eqref{eq:model2}, and let $Y$ and $Z$ be
defined by \eqref{eq:def-Y} and \eqref{eq:def-Z}, respectively. Then,
for every $t>0$,
\begin{equation}\label{eq:SDE-Y-int}
Y_t = X_0t^{\frac12-H}
+ \int_0^t(aZ_s-bY_s)\,ds
+ \sigma C_HW_t,
\end{equation}
where
$C_H = B\!\left(\tfrac32-H,H+\tfrac12\right)$.
\end{theorem}

\begin{proof}
We first transform the noise term. By the definition of $B^H$ and the
stochastic Fubini theorem,
\begin{align*}
\int_0^t(t-s)^{\frac12-H}B_s^H\,ds
&=
\int_0^t\!\!\int_u^t
(t-s)^{\frac12-H}(s-u)^{H-\frac12}\,ds\,dW_u \\
&=
C_H\int_0^t(t-u)\,dW_u
=
C_H\int_0^tW_u\,du,
\end{align*}
where the last equality follows from integration by
parts.

Next, for every locally integrable function $f$, the Fubini theorem gives
\begin{equation}\label{eq:transform-drift}
\int_0^t
(t-s)^{\frac12-H}
\left(\int_0^s f(u)\,du\right)ds
=
\frac{1}{\frac32-H}
\int_0^t
(t-u)^{\frac32-H}f(u)\,du.
\end{equation}
Applying this identity to $f=X^{-1}$ and $f=X$, and applying the same
transformation to \eqref{eq:model2}, we obtain
\begin{align}
\int_0^t(t-s)^{\frac12-H}X_s\,ds
&=
\frac{X_0}{\frac32-H}t^{\frac32-H}
+
\frac{a}{\frac32-H}
\int_0^t\frac{(t-u)^{\frac32-H}}{X_u}\,du
\notag\\
&\quad
-
\frac{b}{\frac32-H}
\int_0^t(t-u)^{\frac32-H}X_u\,du
+
\sigma C_H\int_0^tW_u\,du.
\label{eq:fractional_identity}
\end{align}
Differentiating with respect to $t$ gives
\begin{equation}\label{eq:SDE-Y0}
Y_t
=
X_0t^{\frac12-H}
+
a\int_0^t\frac{(t-u)^{\frac12-H}}{X_u}\,du
-
b\int_0^t(t-u)^{\frac12-H}X_u\,du
+
\sigma C_HW_t.
\end{equation}
On the other hand, integrating
\eqref{eq:def-Y}--\eqref{eq:def-Z} from $0$ to $t$ yields
\begin{equation}\label{eq:YZ-integrals}
\int_0^t(t-u)^{\frac12-H}X_u\,du
=
\int_0^tY_s\,ds,
\qquad
\int_0^t\frac{(t-u)^{\frac12-H}}{X_u}\,du
=
\int_0^tZ_s\,ds.
\end{equation}
Substitution into \eqref{eq:SDE-Y0} proves
\eqref{eq:SDE-Y-int}.
\end{proof}

The following square-integrability property will be used below to
justify the stochastic integrals appearing in the quasi-score
construction. It follows from representations
\eqref{forY}--\eqref{forZ} and the moment bounds of
Remark~\ref{remark4}.

\begin{lemma}\label{lem:1}
Let $X$ be the solution of \eqref{eq:model2}, and let $Y$ and $Z$ be
defined by \eqref{eq:def-Y}--\eqref{eq:def-Z}. Then, for every $T>0$,
\[
\M\int_0^T Y_t^2\,dt<\infty,
\qquad
\M\int_0^T Z_t^2\,dt<\infty.
\]
\end{lemma}

\begin{proof}
Fix $\lambda\in(1/2,H)$. By Remark~\ref{remark4}, there exists a random
variable $\Upsilon$ such that
\[
|X_t-X_s|
\le
\Upsilon|t-s|^\lambda,
\qquad s,t\in[0,T],
\]
and $\M\Upsilon^r<\infty$ for every $r>0$. Moreover,
\eqref{eq:moments} gives finite moments of all orders for
$\sup_{t\in[0,T]}X_t$ and $\sup_{t\in[0,T]}X_t^{-1}$.

We first bound the singular terms in \eqref{forY} and \eqref{forZ}.
We have
\[
\M\int_0^T
\left(
\frac{X_t}{t^{H-\frac12}}
\right)^2dt
\le
\M\sup_{t\in[0,T]}X_t^2
\int_0^Tt^{1-2H}\,dt
=
\frac{T^{2-2H}}{2-2H}
\M\sup_{t\in[0,T]}X_t^2
<\infty.
\]
Similarly,
\[
\M\int_0^T
\left(
\frac{X_t^{-1}}{t^{H-\frac12}}
\right)^2dt
\le
\frac{T^{2-2H}}{2-2H}
\M\sup_{t\in[0,T]}X_t^{-2}
<\infty.
\]

For the integral term in \eqref{forY}, H\"older continuity gives
\[
\left|
\int_0^t
\frac{X_t-X_s}{(t-s)^{H+\frac12}}\,ds
\right|
\le
\Upsilon
\int_0^t
(t-s)^{\lambda-H-\frac12}\,ds
=
\frac{\Upsilon\,t^{\lambda-H+\frac12}}
     {\lambda-H+\frac12}.
\]
Since $\lambda-H+\frac12>0$,
\[
\M\int_0^T
\left(
\int_0^t
\frac{X_t-X_s}{(t-s)^{H+\frac12}}\,ds
\right)^2dt
\le
\frac{\M\Upsilon^2}
     {\left(\lambda-H+\frac12\right)^2}
\frac{T^{2\lambda-2H+2}}
     {2\lambda-2H+2}
<\infty.
\]

For the integral term in \eqref{forZ}, observe that
\[
|X_t^{-1}-X_s^{-1}|
=
\frac{|X_t-X_s|}{X_tX_s}
\le
\Upsilon
\left(
\sup_{u\in[0,T]}X_u^{-2}
\right)
|t-s|^\lambda.
\]
Therefore,
\[
\left|
\int_0^t
\frac{X_t^{-1}-X_s^{-1}}
     {(t-s)^{H+\frac12}}\,ds
\right|
\le
\frac{\Upsilon\,t^{\lambda-H+\frac12}}
     {\lambda-H+\frac12}
\sup_{u\in[0,T]}X_u^{-2}.
\]
Consequently,
\begin{align*}
&\M\int_0^T
\left(
\int_0^t
\frac{X_t^{-1}-X_s^{-1}}
     {(t-s)^{H+\frac12}}\,ds
\right)^2dt
\le
\frac{T^{2\lambda-2H+2}}
     {\left(\lambda-H+\frac12\right)^2
      (2\lambda-2H+2)}
\M\left[
\Upsilon^2
\sup_{u\in[0,T]}X_u^{-4}
\right]
\\
&\qquad\le
\frac{T^{2\lambda-2H+2}}
     {\left(\lambda-H+\frac12\right)^2
      (2\lambda-2H+2)}
\left(\M\Upsilon^4\right)^{1/2}
\left(
\M\sup_{u\in[0,T]}X_u^{-8}
\right)^{1/2}
<\infty,
\end{align*}
where the last inequality follows from the Cauchy--Schwarz inequality
and the moment bounds \eqref{momentgamma} and \eqref{eq:moments}.

Combining these estimates with
\eqref{forY}--\eqref{forZ} proves the assertion.
\end{proof}

\subsection{Quasi-score and estimator}

For the quasi-likelihood construction, it is convenient to remove the
singular deterministic term in \eqref{eq:SDE-Y-int}. Define
\begin{equation}\label{eq:Ybar}
\overline Y_t
\coloneqq
Y_t-X_0t^{\frac12-H},
\qquad t>0,
\end{equation}
and set $\overline Y_0=0$. Since $X_0$ and $H$ are known and $Y$ is
observable, $\overline Y$ is observable as well. By
\eqref{eq:SDE-Y-int},
\begin{equation}\label{eq:Ybar-SDE}
\overline Y_t
=
\int_0^t(aZ_s-bY_s)\,ds
+
\sigma C_HW_t,
\qquad t\in[0,T].
\end{equation}

Introduce
\begin{gather}\label{eq:ABD}
A_T
\coloneqq
\int_0^T Z_t^2\,dt,
\qquad
B_T
\coloneqq
\int_0^T Y_tZ_t\,dt,
\qquad
D_T
\coloneqq
\int_0^T Y_t^2\,dt,
\shortintertext{and}
\label{eq:UV}
U_T
\coloneqq
\int_0^T Z_t\,d\overline Y_t,
\qquad
V_T
\coloneqq
\int_0^T Y_t\,d\overline Y_t.
\end{gather}
We first verify that these quantities are well defined and establish
the nondegeneracy property needed below.

\begin{lemma}\label{lem:estimation-quantities}
For every $T>0$, the random variables
$A_T$, $B_T$, $D_T$, $U_T$, and $V_T$ defined by
\eqref{eq:ABD}--\eqref{eq:UV} are well defined and belong to
$L^1(\mathbb P)$. Moreover,
\begin{equation}\label{eq:Delta-positive}
\Delta_T
\coloneqq
A_TD_T-B_T^2
>0
\qquad\text{a.s.}
\end{equation}
\end{lemma}

\begin{proof}
By Lemma~\ref{lem:1},
\[
\M A_T
=
\M\int_0^T Z_t^2\,dt
<\infty,
\qquad
\M D_T
=
\M\int_0^T Y_t^2\,dt
<\infty.
\]
Thus $A_T$ and $D_T$ are well defined and integrable. By the
Cauchy--Schwarz inequality,
\[
B_T^2\le A_TD_T,
\qquad
|B_T|\le(A_TD_T)^{1/2}.
\]
Hence
\[
\M|B_T|
\le
(\M A_T)^{1/2}(\M D_T)^{1/2}
<\infty,
\]
and therefore $\Delta_T\ge0$ a.s.

Next, by \eqref{eq:Ybar-SDE},
\begin{align}
U_T
&=
aA_T-bB_T
+
\sigma C_H\int_0^T Z_t\,dW_t,
\label{eq:U-rep}\\
V_T
&=
aB_T-bD_T
+
\sigma C_H\int_0^T Y_t\,dW_t.
\label{eq:V-rep}
\end{align}
The processes $Y$ and $Z$ are progressively measurable, and
Lemma~\ref{lem:1} gives
\[
\M\int_0^T Y_t^2\,dt<\infty,
\qquad
\M\int_0^T Z_t^2\,dt<\infty.
\]
Hence the It\^o integrals in
\eqref{eq:U-rep}--\eqref{eq:V-rep} are well defined and square
integrable. The remaining terms are integrable by the preceding
estimates. Thus $U_T$ and $V_T$ are well defined and integrable.

It remains to show that $\Delta_T>0$ a.s. First observe that
$A_T>0$ and $D_T>0$ a.s. Indeed, if $D_T=0$, then $Y_t=0$ for almost
every $t\in[0,T]$, and \eqref{eq:YZ-integrals} would imply
\[
\int_0^t
(t-s)^{\frac12-H}X_s\,ds
=
0,
\qquad t\in[0,T],
\]
which is impossible for $t>0$ because $X_s>0$. The same argument
applied to $Z$ shows that $A_T>0$.

Suppose now that $\Delta_T=0$. Then equality holds in the
Cauchy--Schwarz inequality, so there exists a constant $c=c(\omega)$
such that
\[
Z_t=cY_t
\qquad\text{for almost every }t\in[0,T].
\]
Integrating this relation from $0$ to $t$ and using
\eqref{eq:YZ-integrals}, we obtain
\[
\int_0^t
(t-s)^{\frac12-H}
\left(
\frac1{X_s}-cX_s
\right)ds
=
0,
\qquad t\in[0,T].
\]
By the injectivity of the Riemann--Liouville fractional integral
$I_{0+}^{\frac32-H}$, see \cite[Lemma~2.5]{SMK},
\[
\frac1{X_t}=cX_t
\qquad\text{for almost every }t\in[0,T].
\]
Since $X$ is continuous and strictly positive, this relation holds for
every $t\in[0,T]$. Hence $X$ is constant and necessarily
$X_t\equiv X_0$.

Substitution into \eqref{eq:model2} gives
\[
\sigma B_t^H
=
\left(
bX_0-\frac{a}{X_0}
\right)t,
\qquad t\in[0,T].
\]
In particular, on the event $\{\Delta_T=0\}$,
\[
B_T^H
=
\frac1{\sigma}
\left(
bX_0-\frac{a}{X_0}
\right)T.
\]
Since $B_T^H$ is a nondegenerate Gaussian random variable,
$\mathbb P(\Delta_T=0)=0$. This proves
\eqref{eq:Delta-positive}.
\end{proof}

Let
$\Theta=(0,\infty)^2$,
$\vartheta=(a,b)^\top\in\Theta$,
and define
\[
F_t
\coloneqq
\begin{pmatrix}
Z_t\\
-Y_t
\end{pmatrix},
\qquad t\in[0,T].
\]
Then \eqref{eq:Ybar-SDE} can be written as
\begin{equation}\label{eq:Ybar-vector}
d\overline Y_t
=
F_t^\top\vartheta\,dt
+
\sigma C_H\,dW_t.
\end{equation}
For a fixed observed trajectory of $X$, the process $F$ is completely
determined by the observation and does not depend on the candidate
value of the parameter.

For $\vartheta\in\Theta$, set
\begin{equation}\label{eq:martingale-residual}
M_t(\vartheta)
\coloneqq
\overline Y_t
-
\int_0^tF_s^\top\vartheta\,ds.
\end{equation}
Thus \eqref{eq:Ybar-vector} has the semimartingale decomposition
\[
d\overline Y_t
=
dA_t(\vartheta)+dM_t(\vartheta),
\qquad
dA_t(\vartheta)
=
F_t^\top\vartheta\,dt.
\]
At the true parameter value,
\[
M_t(\vartheta)
=
\sigma C_HW_t,
\qquad
d\langle M(\vartheta)\rangle_t
=
\sigma^2C_H^2\,dt.
\]

We now apply the martingale quasi-likelihood approach developed in
\cite{Heyde1994}; see also \cite[Sect.~8.3]{Heyde1997}. According to
formula~(8.7) in \cite{Heyde1997}, the quasi-score is
\[
S_T(\vartheta)
=
\int_0^T
\M\!\left[
d\dot M_t(\vartheta)\mid\mathcal F_{t-}
\right]
\bigl(d\langle M(\vartheta)\rangle_t\bigr)^{-1}
\,dM_t(\vartheta),
\]
where $\dot M(\vartheta)$ denotes differentiation with respect to the
parameter vector $\vartheta$.

Since $F_t$ does not depend on the candidate value of $\vartheta$,
\[
d\dot M_t(\vartheta)
=
-\frac{\partial}{\partial\vartheta}
\bigl(F_t^\top\vartheta\bigr)\,dt
=
-F_t\,dt.
\]
As $F_t$ is predictable,
\[
\M\!\left[
d\dot M_t(\vartheta)\mid\mathcal F_{t-}
\right]
=
-F_t\,dt.
\]
Therefore,
\[
S_T(\vartheta)
=
-\frac1{\sigma^2C_H^2}
\int_0^T F_t\,dM_t(\vartheta),
\]
and, by \eqref{eq:martingale-residual},
\begin{equation}\label{eq:quasi-score}
S_T(\vartheta)
=
-\frac1{\sigma^2C_H^2}
\int_0^T
F_t
\bigl(
d\overline Y_t-F_t^\top\vartheta\,dt
\bigr).
\end{equation}

By Lemma~\ref{lem:1} and \eqref{eq:Ybar-SDE}, all integrals in
\eqref{eq:quasi-score} are well defined. At the true parameter value,
\[
S_T(\vartheta)
=
-\frac1{\sigma C_H}
\int_0^T F_t\,dW_t.
\]
Hence $\{S_t(\vartheta),\,t\in[0,T]\}$ is a centered
square-integrable martingale and, in particular,
$\M S_T(\vartheta)=0$.
Using \eqref{eq:ABD}--\eqref{eq:UV}, the quasi-score can be written as
\begin{equation}\label{eq:quasi-score-components}
S_T(a,b)
=
-\frac1{\sigma^2C_H^2}
\begin{pmatrix}
U_T-aA_T+bB_T\\[1mm]
-V_T+aB_T-bD_T
\end{pmatrix}.
\end{equation}

The quasi-likelihood estimator is defined as a solution of the
estimating equation
\begin{equation}\label{eq:quasi-score-equation}
S_T(\vartheta)=0.
\end{equation}

\begin{theorem}\label{thm:QLE}
For every $T>0$, equation \eqref{eq:quasi-score-equation}, considered
on $\mathbb R^2$, has a.s.\ a unique solution
$(\widehat a_T,\widehat b_T)$ given by
\begin{equation}\label{eq:QLE-ab-final}
\widehat a_T
=
\frac{D_TU_T-B_TV_T}{\Delta_T},
\qquad
\widehat b_T
=
\frac{B_TU_T-A_TV_T}{\Delta_T},
\end{equation}
where
\[
\Delta_T=A_TD_T-B_T^2.
\]
Consequently, if
$(\widehat a_T,\widehat b_T)\in\Theta$, then it is the unique
quasi-likelihood estimator of $(a,b)$ in the parameter space $\Theta$.
\end{theorem}

\begin{proof}
By \eqref{eq:quasi-score-components},
$S_T(a,b)=0$ is equivalent to
\[
\begin{cases}
A_Ta-B_Tb=U_T,\\
B_Ta-D_Tb=V_T.
\end{cases}
\]
The determinant of this system equals
$A_TD_T-B_T^2 = \Delta_T$.
By Lemma~\ref{lem:estimation-quantities},
$\Delta_T>0$ a.s., so the system has a.s.\ a unique solution.
Solving it yields \eqref{eq:QLE-ab-final}.
\end{proof}

\begin{remark}[Quasi-likelihood function]
\label{rem:quasilik-function}
Although the quasi-likelihood estimator is defined directly through the
estimating equation $S_T(\vartheta)=0$, in the present model the
quasi-score can also be associated with a scalar quasi-likelihood
function. Define
\begin{equation}\label{eq:quasilik}
Q_T(a,b)
=
\frac1{\sigma^2C_H^2}
\left(
aU_T-bV_T
-\frac12a^2A_T
+abB_T
-\frac12b^2D_T
\right).
\end{equation}
Differentiating with respect to $a$ and $b$ and comparing with
\eqref{eq:quasi-score-components}, we obtain
\[
\nabla Q_T(\vartheta)
=
-S_T(\vartheta).
\]
Moreover,
\[
\nabla^2Q_T
=
-\frac1{\sigma^2C_H^2}
\begin{pmatrix}
A_T&-B_T\\
-B_T&D_T
\end{pmatrix}.
\]
By Lemma~\ref{lem:estimation-quantities},
$\nabla^2Q_T$ is negative definite a.s. Hence $Q_T$ is strictly
concave, and the solution \eqref{eq:QLE-ab-final} of the quasi-score
equation is equivalently the unique maximizer of $Q_T$ over
$\mathbb R^2$.
\end{remark}

\begin{remark}[Relation to maximum likelihood estimation]
\label{rem:QLE-MLE}
The function $Q_T$ also clarifies the relation between quasi-likelihood
and maximum likelihood estimation; see \cite{Heyde1994} and
\cite[Sect.~8.3]{Heyde1997} for a general discussion. Suppose that,
for $\vartheta,\vartheta_0\in\Theta$, the laws
$\mathbb P_\vartheta^T$ and $\mathbb P_{\vartheta_0}^T$ of the
transformed observation $\overline Y$ on $[0,T]$ are mutually
absolutely continuous and that the usual Girsanov likelihood formula is
applicable. Then
\[
\log
\frac{d\mathbb P_\vartheta^T}
     {d\mathbb P_{\vartheta_0}^T}
=
Q_T(\vartheta)-Q_T(\vartheta_0).
\]
Thus, whenever the ordinary likelihood is well defined, the
quasi-likelihood estimator \eqref{eq:QLE-ab-final} coincides with the
corresponding maximum likelihood estimator based on $\overline Y$.
In contrast, the quasi-likelihood construction itself does not require
mutual absolute continuity of the parameterized family of laws or any
additional exponential-martingale condition.
\end{remark}

\begin{remark}[Minimum contrast interpretation]
The estimator \eqref{eq:QLE-ab-final} also admits a minimum contrast
interpretation; see, for example, \cite{Lanska1979} for minimum
contrast estimation in diffusion models. Define the quadratic contrast
criterion
\[
\mathcal C_T(a,b)
=
a^2A_T-2abB_T+b^2D_T-2aU_T+2bV_T.
\]
Then
$\mathcal C_T(a,b) = -2\sigma^2C_H^2\,Q_T(a,b)$,
and hence minimizing $\mathcal C_T$ is equivalent to maximizing
$Q_T$.
\end{remark}

\section{Asymptotic estimation of the drift-parameter ratio}
\label{sec:ratio}

In Section~\ref{sec2}, we constructed explicit quasi-likelihood
estimators of both drift parameters in the case of
Riemann--Liouville fractional Brownian noise. We now consider a
complementary estimation problem that can be treated under substantially
more general assumptions on the Gaussian driving process.

More precisely, we study the ratio
$\theta:=\frac{a}{b}$,
which has a natural interpretation in connection with the
mean-reverting CIR-type structure of the model discussed in the
Introduction. We show that, under suitable regularity and growth
conditions on the driving process,
\[
\widehat\theta_T
:=
\frac{\int_0^T X_s\,ds}
     {\int_0^T X_s^{-1}\,ds}
\longrightarrow
\theta
\qquad\text{a.s. as }T\to\infty.
\]
Thus $\widehat\theta_T$ provides a strongly consistent estimator of the
parameter ratio $a/b$. Unlike the quasi-likelihood construction of
Section~\ref{sec2}, this estimator does not rely on a specific Volterra
representation of the noise and applies to a broad class of Gaussian
driving processes. We verify the required assumptions for several
important examples below. The finite-sample behavior of the estimator
is investigated numerically in Subsections~\ref{sec:sim-ratio-design} and \ref{sec:sim-ratio-results}.

\subsection{Strong consistency of the ratio estimator}
\label{sec:relation}

To establish the asymptotic behavior of $\widehat\theta_T$, we first
consider the auxiliary Bessel-type equation
\begin{equation}\label{bessel}
\widetilde X_t
=
X_0
+
a\int_0^t\frac{ds}{\widetilde X_s}
+
G_t,
\qquad t\ge0,
\end{equation}
where $X_0>0$, $a>0$, and $G$ is the same driving process as in
\eqref{eq:model}, with $G_0=0$. Equation \eqref{bessel} is the
particular case of \eqref{eq:model} corresponding to $b=0$.
The following result is an immediate consequence of
Theorem~\ref{theorem1}.

\begin{theorem}\label{theorem2}
Let $T>0$, and suppose that $G$ is $\lambda$-H\"older continuous on
$[0,T]$ for some $\lambda>1/2$. Then equation \eqref{bessel} has a
unique strictly positive solution on $[0,T]$.
\end{theorem}

\begin{remark}\label{rem:comparison}
Suppose that the assumptions of Theorem~\ref{theorem2} hold. Let $X$
and $\widetilde X$ denote the solutions of \eqref{eq:model} and
\eqref{bessel}, respectively, with the same initial value $X_0$, the
same coefficient $a$, and driven by the same process $G$. Then
\[
0<X_t\le\widetilde X_t,
\qquad t\in[0,T],
\quad\text{a.s.}
\]
This follows from the general comparison result
\cite[Lemma~3.2]{MishYur3}; see also \cite[Section~4]{MishYur1} for
the corresponding result in the fractional Brownian setting.
\end{remark}

For a function $g\colon[0,\infty)\to\R$, denote
\[
g_t^*
:=
\sup_{0\le s\le t}|g_s|,
\qquad t\ge0.
\]

We now establish the strong consistency of the ratio estimator under a
sublinear-growth assumption on the driving process.

\begin{theorem}\label{thm:parameter-relation}
Let $G$ be the driving process in \eqref{eq:model}. Assume that $G$ is
$\lambda$-H\"older continuous on every interval $[0,T]$ for some
$\lambda>1/2$ and that
\begin{equation}\label{eq:G-growth}
\frac{G_t^*}{t}
\to0
\qquad\text{a.s. as }t\to\infty.
\end{equation}
Then
\begin{equation}\label{eq:theta-consistency}
\widehat\theta_T
\coloneqq
\frac{\int_0^T X_s\,ds}
     {\int_0^T X_s^{-1}\,ds}
\to
\frac{a}{b}
\qquad\text{a.s. as }T\to\infty.
\end{equation}
Thus $\widehat\theta_T$ is a strongly consistent estimator of
$\theta=a/b$.
\end{theorem}

\begin{proof}
We first establish a sublinear growth bound for the auxiliary process $\widetilde X$. 
We work on an event of probability one on which \eqref{eq:G-growth} holds and $G$ is $\lambda$-H\"older continuous on every interval $[0,T]$.
Arguing as in Subsection~\ref{subsec:transformation}, one can show that strict positivity of
$\widetilde X$ and \eqref{bessel} imply that $\widetilde X$ is also
$\lambda$-H\"older continuous on every interval $[0,T]$.

Since $\lambda>1/2$, the integrals
$\int_0^t \widetilde X_s\,dG_s$,
$\int_0^t G_s\,d\widetilde X_s$, and
$\int_0^t G_s\,dG_s$
are well defined in the Young sense; see \cite{Young1936} and
\cite[Sect.~4.1]{FrizHairer}.

By the change-of-variable and integration-by-parts formulas for Young
integrals, we obtain
\begin{align*}
\widetilde X_t^2
&= X_0^2+2\int_0^t \widetilde X_s\,d\widetilde X_s
 = X_0^2+2at+2\int_0^t\widetilde X_s\,dG_s 
= X_0^2+2at+2\widetilde X_tG_t
   -2\int_0^tG_s\,d\widetilde X_s.
\end{align*}
Therefore, using \eqref{bessel}, we get the identities
\begin{align*}
\widetilde X_t^2
&= X_0^2+2at+2\widetilde X_tG_t
-2a\int_0^t\frac{G_s}{\widetilde X_s}\,ds
-2\int_0^tG_s\,dG_s \\
&= X_0^2+2at+2\widetilde X_tG_t
-2a\int_0^t\frac{G_s}{\widetilde X_s}\,ds
-G_t^2.
\end{align*}

Since $|G_s|\le G_t^*$ for $s\le t$, we obtain that
\[
\widetilde X_t^2
\le
X_0^2+2at+2\widetilde X_tG_t^*
+2aG_t^*\int_0^t\frac{ds}{\widetilde X_s}.
\]
By \eqref{bessel},
\[
a\int_0^t\frac{ds}{\widetilde X_s}
=
\widetilde X_t-X_0-G_t
\le \widetilde X_t+G_t^*.
\]
Consequently,
\[
\widetilde X_t^2
\le
X_0^2+2at+4G_t^*\widetilde X_t+2(G_t^*)^2.
\]

Thus,
\[
\widetilde X_t^2
-4G_t^*\widetilde X_t
-\bigl(X_0^2+2at+2(G_t^*)^2\bigr)
\le0.
\]
Since $\widetilde X_t>0$, solving this quadratic inequality yields
\[
\widetilde X_t
\le
2G_t^*
+\sqrt{X_0^2+2at+6(G_t^*)^2}.
\]
Therefore,
\[
0\le \frac{\widetilde X_t}{t}
\le
2\frac{G_t^*}{t}
+
\sqrt{
\frac{X_0^2}{t^2}
+\frac{2a}{t}
+6\left(\frac{G_t^*}{t}\right)^2
}.
\]
By \eqref{eq:G-growth},
$\frac{G_t^*}{t}\to0$
a.s. as $t\to\infty$.
Hence
\[
\frac{\widetilde X_t}{t}\to0
\quad\text{a.s. as }t\to\infty.
\]

Applying Remark~\ref{rem:comparison} on every finite interval, we have
$0<X_t\le\widetilde X_t$, $t\ge0$, a.s.
Hence
$\frac{X_t}{t}\to0$ a.s.\ as $t\to\infty$.
Therefore, dividing \eqref{eq:model} by $t$, we obtain
\[
R_t\coloneqq
\frac1t\left(
a\int_0^t\frac{ds}{X_s}
-b\int_0^t X_s\,ds
\right)
=
\frac{X_t-X_0-G_t}{t}
\to0
\quad\text{a.s.\ as }t\to\infty.
\]

By the Cauchy--Schwarz inequality,
\[
\left(\int_0^t\frac{ds}{X_s}\right)
\left(\int_0^t X_s\,ds\right)\ge t^2.
\]
Consequently,
\[
a\left(\frac1t\int_0^t\frac{ds}{X_s}\right)^2
-R_t\left(\frac1t\int_0^t\frac{ds}{X_s}\right)
-b\ge0.
\]
Since $\frac1t\int_0^t X_s^{-1}\,ds>0$, solving this quadratic inequality gives
\[
\frac1t\int_0^t\frac{ds}{X_s}
\ge
\frac{R_t+\sqrt{R_t^2+4ab}}{2a}.
\]
Therefore,
\[
\liminf_{t\to\infty}
\frac1t\int_0^t\frac{ds}{X_s}
\ge\sqrt{\frac ba}>0
\quad\text{a.s.}
\]
Finally,
\[
a-b\,
\frac{\int_0^t X_s\,ds}
     {\int_0^t X_s^{-1}\,ds}
=
\frac{R_t}
{\frac1t\int_0^t X_s^{-1}\,ds}
\to0
\quad\text{a.s. as }t\to\infty,
\]
which proves the claim.
\end{proof}

\subsection{Examples of admissible Gaussian driving processes}
\label{subsec:ratio-examples}

The literature contains numerous results on the asymptotic growth of
Gaussian processes under various structural assumptions. In particular,
such estimates are available for stationary processes
\cite[Section~1]{Marcus72}, processes with stationary increments
\cite[Section~2]{Marcus72}, self-similar processes
\cite[Theorem~2.1]{MRS2023}, and several important particular classes.
Together with suitable regularity estimates, these results allow us to
verify the assumptions of Theorem~\ref{thm:parameter-relation} for a
broad range of Gaussian driving noises.

\begin{example}[Fractional Brownian motion]
Let $G=W^H$ be a fractional Brownian motion with Hurst parameter
$H\in(1/2,1)$. It is well known that, for every $\lambda\in(0,H)$,
the sample paths of $W^H$ are locally $\lambda$-H\"older continuous
almost surely; see, e.g., \cite{Mishura2008}. Hence the H\"older
continuity assumption of Theorem~\ref{thm:parameter-relation} is
satisfied by choosing any $\lambda\in(1/2,H)$.

Furthermore, according to
\cite[Remark~7]{KozachenkoMelnikovMishura2015}, for every $p>1$,
\[
(W^H)^*_t
\coloneqq
\sup_{0\le s\le t}|W_s^H|
\le
\bigl(t^H(\log t)^p\vee1\bigr)\,\xi(p)
\]
for all sufficiently large $t$, where $\xi(p)$ is an a.s.\ finite
random variable independent of $t$. Since $H<1$, it follows that
\[
\frac{(W^H)^*_t}{t}
\le
\xi(p)\,t^{H-1}(\log t)^p
\to0
\qquad\text{a.s.}
\]
Thus condition \eqref{eq:G-growth} is satisfied, and
Theorem~\ref{thm:parameter-relation} can be applied.
\end{example}

The preceding growth estimate for fractional Brownian motion is a
particular instance of a more general result for self-similar Gaussian
processes. In particular, \cite[Theorem~2.1]{MRS2023} provides an
asymptotic growth bound for any continuous centered self-similar
Gaussian process. This allows us to treat a broader class of Gaussian
Volterra processes in the same way.

\begin{example}[Gaussian Volterra process with power-type kernel]
Consider the Gaussian Volterra process introduced in
\cite{MishuraShklyar2022I},
\[
G_t
=
\int_0^t
s^\alpha
\left(
\int_s^t u^\beta (u-s)^\gamma\,du
\right)dW_s,
\qquad t\ge0,
\]
where
\[
\alpha>-\tfrac12,\qquad
\gamma>-1,\qquad
\alpha+\beta+\gamma>-\tfrac32.
\]
Under these conditions the process is well defined; see
\cite[Theorem~1]{MishuraShklyar2022I}. Moreover, it is self-similar
with exponent
\[
h=\alpha+\beta+\gamma+\tfrac32;
\]
see \cite[Proposition~1]{MishuraShklyar2022I}.

For the purposes of Theorem~\ref{thm:parameter-relation}, assume
additionally that
\[
-1<\alpha+\beta+\gamma<-\tfrac12,
\]
or, equivalently,
\[
\tfrac12<h<1.
\]
According to \cite[Theorem~6]{MishuraShklyar2022I}, on every compact
interval $[0,T]$ the process admits a modification whose paths are
H\"older continuous of every order
\[
\lambda < \min\left\{h,\,\gamma+\tfrac32,\,1\right\}.
\]
Since $h>1/2$ and $\gamma>-1$, the quantity on the right-hand side is
strictly greater than $1/2$. Hence the H\"older continuity assumption
of Theorem~\ref{thm:parameter-relation} is satisfied.

Furthermore, \cite[Theorem~2.3]{MRS2023} yields
\[
\sup_{t\ge3}
\frac{|G_t|}
     {t^h\sqrt{\log\log t}}
<\infty
\qquad\text{a.s.}
\]
Since $h>0$, this implies
\[
G_t^*
=
O\left(t^h\sqrt{\log\log t}\right)
\qquad\text{a.s.}
\]
Therefore, since $h<1$,
\[
\frac{G_t^*}{t}
=
O\left(t^{h-1}\sqrt{\log\log t}\right)
\to0
\qquad\text{a.s.}
\]
Thus all assumptions of Theorem~\ref{thm:parameter-relation} are
satisfied.

Both fractional Brownian motion and Riemann--Liouville fractional
Brownian motion arise as particular cases of this class. Indeed, for
$H\in(1/2,1)$, fractional Brownian motion is obtained, up to a
deterministic normalizing constant, by taking
\[
\alpha=\tfrac12-H,\qquad
\beta=H-\tfrac12,\qquad
\gamma=H-\tfrac32.
\]
This choice gives the classical Volterra representation of fractional
Brownian motion; see \cite{MishuraShklyar2022I}.

On the other hand, taking
\[
\alpha=\beta=0,\qquad
\gamma=H-\tfrac32,
\]
we obtain
\[
\int_0^t
\left(
\int_s^t (u-s)^{H-\frac32}\,du
\right)dW_s
=
\frac{1}{H-\frac12}
\int_0^t(t-s)^{H-\frac12}\,dW_s.
\]
Thus, up to the factor $(H-\frac12)^{-1}$, this is the
Riemann--Liouville fractional Brownian motion with Hurst parameter $H$.
In particular, both the fractional Brownian motion considered above
and the Riemann--Liouville fractional Brownian motion used in
Section~\ref{sec2} belong to this broader class.
\end{example}

\begin{example}[Tempered fractional Brownian motion]
Let $G=B^{\mathrm{I}}_{H,\lambda}$ be a tempered fractional Brownian
motion with parameters $H>1/2$ and $\lambda>0$, introduced in
\cite{MeerschaertSabzikar2013}. It is defined by
\[
B^{\mathrm{I}}_{H,\lambda}(t)
=
\int_{\mathbb R}
\left(
e^{-\lambda(t-s)_+}(t-s)_+^{H-\frac12}
-
e^{-\lambda(-s)_+}(-s)_+^{H-\frac12}
\right)dW_s,
\qquad t\in\mathbb R,
\]
where $x_+=\max\{x,0\}$ and $W$ is a two-sided Wiener process.

If $H\in(0,1)$, the sample paths of
$B^{\mathrm{I}}_{H,\lambda}$ are locally H\"older continuous of every
order $\eta<H$; see
\cite[Theorem~5.1]{MeerschaertSabzikar2013}.
For $H=1$, they are locally H\"older continuous of every order
$\eta<1$, while for $H>1$ they are continuously differentiable;
see also \cite[Remark~3]{MishuraRalchenko2024}.
Thus, for every $H>1/2$, the H\"older continuity assumption of
Theorem~\ref{thm:parameter-relation} is satisfied.

According to \cite[Theorem~3]{MishuraRalchenko2024}, for every
$\delta>0$ there exists an a.s.\ finite random variable
$\xi=\xi(\delta)$ such that
\[
\sup_{0\le s\le t}
\left|B^{\mathrm{I}}_{H,\lambda}(s)\right|
\le
(t^\delta\vee1)\xi,
\qquad t>0,
\quad\text{a.s.}
\]
Choosing $\delta\in(0,1)$, we obtain
\[
\frac{\bigl(B^{\mathrm{I}}_{H,\lambda}\bigr)^*_t}{t}
\le
\xi\,t^{\delta-1}
\to0
\qquad\text{a.s.}
\]
Hence condition \eqref{eq:G-growth} is satisfied, and
Theorem~\ref{thm:parameter-relation} can be applied.
\end{example}

\begin{example}[Tempered fractional Brownian motion of the second kind]
Let $G=B^{\mathrm{II}}_{H,\lambda}$ be a tempered fractional Brownian
motion of the second kind with $H>1/2$ and $\lambda>0$, introduced in
\cite{SabzikarSurgailis2018}. It is defined by
\[
B^{\mathrm{II}}_{H,\lambda}(t)
=
\int_{\mathbb R}
g^{\mathrm{II}}_{H,\lambda,t}(s)\,dW_s,
\qquad t\in\mathbb R,
\]
where
\[
g^{\mathrm{II}}_{H,\lambda,t}(s)
= e^{-\lambda(t-s)_+}(t-s)_+^{H-\frac12}
- e^{-\lambda(-s)_+}(-s)_+^{H-\frac12}
+ \lambda\int_0^t
(u-s)_+^{H-\frac12}e^{-\lambda(u-s)_+}\,du .
\]
Its sample paths have the same relevant local regularity as those of
TFBM; see \cite[Remark~5]{MishuraRalchenko2024}. In particular, the
H\"older continuity assumption of Theorem~\ref{thm:parameter-relation}
is satisfied whenever $H>1/2$.

The asymptotic growth is different from that of TFBM. According to
\cite[Theorem~5]{MishuraRalchenko2024}, for every $p>1$ there exists
an a.s.\ finite random variable $\xi=\xi(p)$ such that
\[
\sup_{0\le s\le t}
\left|B^{\mathrm{II}}_{H,\lambda}(s)\right|
\le
\left(
t^{1/2}(\log^+ t)^p\vee1
\right)\xi,
\qquad t>0,
\quad\text{a.s.},
\]
where $\log^+t=\max\{\log t,0\}$. Consequently,
\[
\frac{\bigl(B^{\mathrm{II}}_{H,\lambda}\bigr)^*_t}{t}
=
O\left(t^{-1/2}(\log t)^p\right)
\to0
\qquad\text{a.s.}
\]
Thus condition \eqref{eq:G-growth} is satisfied, and
Theorem~\ref{thm:parameter-relation} can be applied in this case as well.
\end{example}

\begin{example}[Multifractional Brownian motion]
Multifractional Brownian motion generalizes fractional Brownian motion
by replacing the constant Hurst parameter with a time-varying Hurst
function. Let
$H\colon\mathbb R_+\to[h_1,h_2]\subset(0,1)$
be a $\kappa$-H\"older continuous function, that is,
\[
|H(t)-H(s)|\le D|t-s|^\kappa,
\qquad s,t\ge0,
\]
for some $D>0$ and $\kappa\in(0,1]$.
The harmonizable multifractional Brownian motion with Hurst function
$H$, introduced in \cite{BenassiJaffardRoux1997}, is defined by
\[
G_t
=
\int_{\mathbb R}
\frac{e^{itu}-1}{|u|^{H(t)+1/2}}\,
\widetilde W(du),
\qquad t\ge0,
\]
where $\widetilde W$ is the Fourier transform of white noise.

Assume that
\[
\min\{h_1,\kappa\}>\frac12.
\]
Then, according to \cite[Remark~3.2]{DKMR18}, the sample paths of $G$
are locally H\"older continuous of every order
$\eta<\min\{h_1,\kappa\}$.
Hence the H\"older continuity assumption of
Theorem~\ref{thm:parameter-relation} is satisfied.

Let
\[
h_*=\limsup_{t\to\infty}H(t).
\]
By \cite[Theorem~3.8]{DKMR18}, for every $\delta>0$ there exist
$t_\delta>0$ and an a.s.\ finite random variable $\xi=\xi(\delta)$
such that
\[
|G_t|
\le
\bigl(t^{h_*+\delta}\vee1\bigr)\xi,
\qquad t>t_\delta,
\quad\text{a.s.}
\]
Since $h_*\le h_2<1$, we can choose
$\delta\in(0,1-h_*)$. For such $\delta$,
$G_t^* = O(t^{h_*+\delta})$ a.s.,
and hence \enlargethispage{10pt}
\[
\frac{G_t^*}{t}
= O\bigl(t^{h_*+\delta-1}\bigr)
\to0
\qquad\text{a.s.}
\]
Thus condition \eqref{eq:G-growth} is satisfied, and
Theorem~\ref{thm:parameter-relation} can be applied.
\end{example}

\section{Simulation study}
\label{sec4}

In this section, we present results of the estimators' performance in numerical simulations.

Subsections~\ref{sec:sim-ab-design} and~\ref{sec:sim-ab-results} investigate the behavior of estimator $(\hat{a}_T, \hat{b}_T)$ of Section~\ref{sec2} under Riemann--Liouville fractional Brownian noise. Subsections~\ref{sec:sim-ratio-design} and \ref{sec:sim-ratio-results} study the behaviour of the ratio estimator $\hat{\theta}_T$ from Section~\ref{sec:ratio} under Riemann--Liouville and standard fractional Brownian noise. All computations were carried out in \texttt{R}; the computationally intensive parts were written in  \texttt{C++} via the ``\texttt{Rcpp}'' library in \texttt{R}.

\subsection{Quasi-likelihood estimator: simulation design}\label{sec:sim-ab-design}

Our model \eqref{eq:model2} is simulated with $a=0.8$, $b=1.2$, $\sigma=0.5$, $X_0=1.5$ and $H\in\{0.6,0.7,0.8,0.9\}$ on the uniform grid
\begin{equation*}
    t_k=k\delta,\qquad k=0,1,\dots,N,\qquad T=N\delta .
\end{equation*}
One Monte Carlo replication consists of four steps: simulation of $B^H$ on the grid, construction of $X$, approximation of the fractional derivatives $Y$ and $Z$, and evaluation of the estimator.

\textbf{Simulation of the driving fractional Brownian motion.} The process $B^H$ is simulated exactly at the grid points from its covariance function (Definition~\ref{def:rlfbm}) using the Cholesky method (see, e.g., \cite[Appendix~3]{BannaMishuraRalchenkoShklyar2019}). 
As a check, the empirical mean squared increments of the simulated trajectories agree with their exact values within about $1\%$.

Exact simulation of $B^H$ by the Cholesky factorization is expensive for large $N$. 
Cheaper schemes with $O(N\log N)$ operations discretize the Volterra representation directly, e.g., by the left-point Riemann sum
$B^H_{t_k}\approx\sum_{j<k}(t_k-t_j)^{H-\frac{1}{2}}(W_{t_{j+1}}-W_{t_j})$ or by a hybrid-type modification that evaluates the kernel at shifted points away from the diagonal (see, e.g.,
\cite[Appendix~3]{BannaMishuraRalchenkoShklyar2019}). 
In our preliminary experiments, both constructions overstated the variance of the increments of $B^H$. 

\textbf{Construction of the process $X$.}
Given the increments $\Delta B_k^H=B^H_{t_{k+1}}-B^H_{t_k}$, we apply the drift-implicit (backward) Euler scheme \cite[Section~4, Example~4.1]{DiNunnoMishuraYT2023} to 
\begin{equation*}
    dX_t=\left(\frac{a}{X_t}-bX_t\right)dt+\sigma\,dB_t^H
\end{equation*}
which produces $X_{k+1}=X_k+\left(\frac{a}{X_{k+1}}-bX_{k+1}\right)\delta +\sigma\Delta B_k^H$. Defining $A=1+b\delta$, $C=a \delta $, $B_k=X_k+\sigma\Delta B_k^H$, the previous equation may be rewritten as
\begin{equation*}
    AX_{k+1}^2-B_kX_{k+1}-C=0.
\end{equation*}
Selecting the positive root yields the explicit recursion
\begin{equation*}
    X_{k+1} =
    \frac{B_k+\sqrt{B_k^2+4AC}}{2A}.
\end{equation*}

The implicit treatment of the drift is crucial for preserving the positivity of the numerical solution.

\textbf{Fractional derivative.}
The estimator is built from $Y=\Gamma(\frac{3}{2}-H)D_{0+}^{\alpha}X$ and $Z=\Gamma(\frac{3}{2}-H)D_{0+}^{\alpha}X^{-1}$, $\alpha=H-\frac{1}{2}$; see \eqref{eq:def-Y}--\eqref{eq:def-Z}.
We approximate the Riemann--Liouville derivative by Lubich's convolution quadrature generated by the second-order backward differentiation formula
(BDF2) \cite[Section~3.2.3]{CaiLi2020}
\begin{equation*}
    D_\delta^{\alpha} f(t_j) = 
    \delta^{-\alpha}\sum_{k=0}^{j} w_k f(t_{j-k}),
    \qquad 
    \sum_{k=0}^{\infty} w_kz^k = 
    \Bigl(\tfrac32 - 2z + \tfrac12z^2\Bigr)^{\alpha}.
\end{equation*}
Differentiating the generating function gives the recursion
\begin{equation*}
    w_0 = \Bigl(\tfrac32\Bigr)^{\alpha},
    \qquad
    w_k = \frac{2}{3k}\Bigl(-2\,(\alpha+1-k)\, 
    w_{k-1}+\tfrac12\,
    (2\alpha+2-k)\,w_{k-2}\Bigr),
    \quad 
    k\ge1,
\end{equation*}
with $w_{-1}=0$. The convolutions are evaluated for all $j$ at once by the fast Fourier transform. Since the Riemann--Liouville derivative of a constant is singular at the origin, $\Gamma(1-\alpha)D_{0+}^{\alpha}c=c\,t^{-\alpha}$, only the zero-start parts are approximated numerically:
   \begin{align*}
    Y_{t_j} & \approx X_0t_j^{-\alpha} + \Gamma(1-\alpha)\, 
               D_\delta^{\alpha}(X-X_0)(t_j),\\
    Z_{t_j} & \approx X_0^{-1}t_j^{-\alpha} + \Gamma(1-\alpha)\,
               D_\delta^{\alpha}\bigl(X^{-1} - X_0^{-1}\bigr)(t_j),
    \qquad j\ge1.
   \end{align*}
Note that $\Gamma(1-\alpha)=\Gamma(\frac{3}{2}-H)$. The approximation is non-anticipating (it uses no future values of the trajectory): $Y_{t_j}$ and $Z_{t_j}$ depend only on $X_{t_0},\dots,X_{t_j}$, which is why we selected it. We also considered the weighted and shifted Gr\"unwald difference (WSGD) operator with shifts $(p,q)=(1,0)$~\cite[Section~3.2.2]{CaiLi2020}, both in its original form and with the analytic zero-start correction used above. 

\textbf{Estimator.} 
Denote $Y_j=Y_{t_j}$, $Z_j=Z_{t_j}$ and $\varphi(t)=X_0\bigl(\tfrac12-H\bigr)t^{-\frac{1}{2}-H}$, the derivative of the deterministic term $X_0t^{\frac{1}{2}-H}$ in \eqref{eq:Ybar}. The quantities \eqref{eq:ABD}--\eqref{eq:UV} are approximated by
\begin{gather*}
    A_T \approx \sum_{j}Z_j^2 \, \delta, \qquad
    B_T \approx \sum_{j}Y_jZ_j\, \delta, \qquad
    D_T \approx \sum_{j}Y_j^2 \, \delta,
\\
U_T \approx \sum_{j}Z_j \bigl(Y_{j+1} -Y_j -\varphi(t_j)\delta\bigr),\qquad
V_T \approx \sum_{j}Y_j \bigl(Y_{j+1} -Y_j -\varphi(t_j)\delta\bigr),
\end{gather*}
and $\widehat a_T$, $\widehat b_T$ are then computed from \eqref{eq:QLE-ab-final}. Since $Y$ and $Z$ retain the singularity $t^{-\alpha}$ at the origin, all sums are restricted to $t_{\mathrm{burn}}\le t_j\le t_{N-2}$ with a fixed burn-in time $t_{\mathrm{burn}}=0.5$.

\textbf{Monte Carlo design.}
We take $a=0.8$, $b=1.2$, $\sigma=0.5$, $X_0=1.5$ and $H\in\{0.6,0.7,0.8$, $0.9\}$. Using the constructions above, we conduct two studies, each with $m=1000$ replications for every $H$. 
In the \emph{step-size study}, we fix the horizon at $T=100$ and take 
$\delta \in \{0.005,$ $0.01,$ $0.02,$ $0.025,$ $0.04,$ $0.05,$ $0.1\}$ 
to evaluate the impact of decreasing the discretization step. 
In the \emph{horizon study}, $\delta \in \{0.01,$ $0.02,$ $0.05,$ $0.1\}$, the estimates are recorded at every integer horizon $T = 5, 6, \dots, 200$, with a focus more on how $T$ affects the accuracy of the estimators. 
In both studies, all step sizes within one replication are obtained from one trajectory of the noise on the finest grid ($N = 2 \cdot 10^4$ points), and all horizons from one trajectory on the longest interval. The sums in the estimators are cumulative, so stopping them at $T$ is the same as recomputing the estimators on $[0,T]$.

\textbf{Accuracy measures.}
For a fixed set of parameters $(H,T,\delta)$, let $\widehat a_T^{(i)}$, $i=1,\dots,m$, be the sample of corresponding estimates obtained in $m$ independent replications. Denote $\overline a_T = \frac1m\sum_{i=1}^m \widehat a_T^{(i)}$ and $s_a^2 = \frac1{m-1} \sum_{i=1}^m\bigl(\widehat a_T^{(i)}-\overline a_T\bigr)^2$.
For our simulations, we evaluate the bias, the standard error, and the root mean squared error:
\begin{equation*}
    \mathrm{Bias}(\widehat a_T) = \overline a_T-a,    \quad
    \mathrm{SE}(\widehat a_T)   = \frac{s_a}{\sqrt m},\quad
    \mathrm{RMSE}(\widehat a_T) = \biggl(\frac1m \sum_{i=1}^m
    \bigl(\widehat a_T^{(i)} - a\bigr)^2\biggr)^{1/2},
\end{equation*}
and the same for the estimate $\widehat b_T$.

\subsection{Quasi-likelihood estimator: results}\label{sec:sim-ab-results}

\begin{figure}%[ht]
\centering
\includegraphics[width=\textwidth]{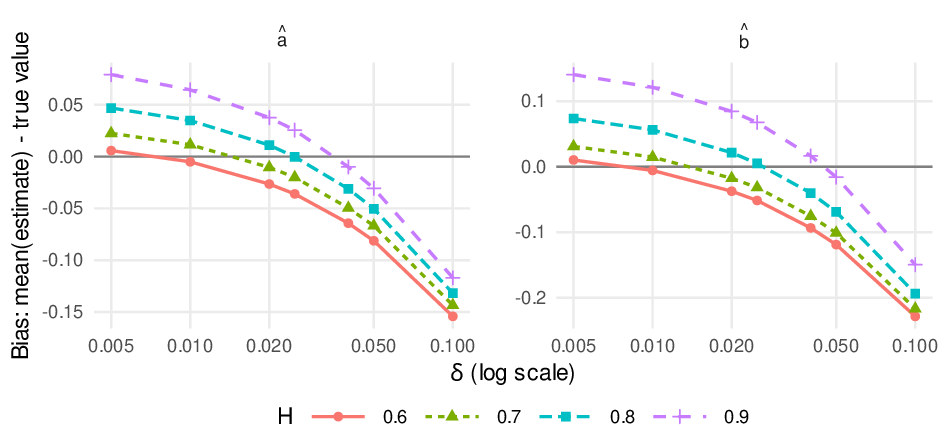}
\caption{Step-size study, $T=100$: bias of $\widehat a_T$ (left) and $\widehat b_T$ (right) as a function of $\delta$ (logarithmic scale), by $H$}
\label{fig:sim-bias-delta}
\end{figure}

\begin{figure}%[ht]
\centering
\includegraphics[width=\textwidth]{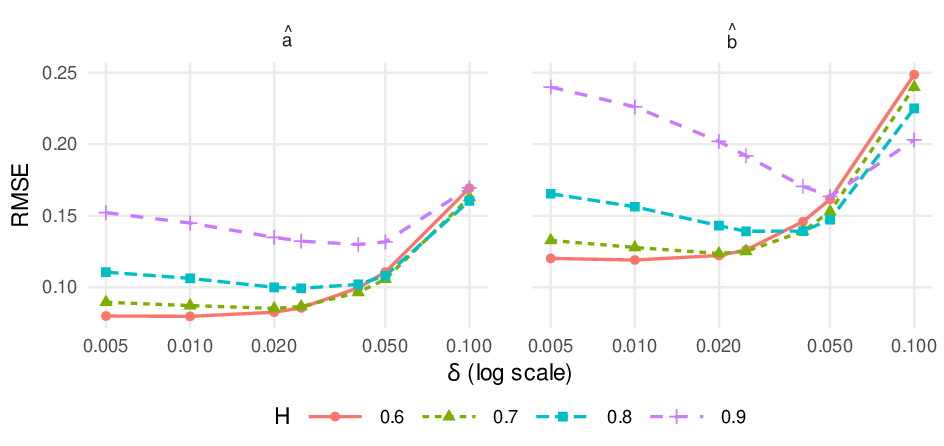}
\caption{Step-size study, $T=100$: RMSE of $\widehat a_T$ (left) and $\widehat b_T$ (right) as a function of $\delta$ (logarithmic scale), by $H$}
\label{fig:sim-rmse-delta}
\end{figure}

\textbf{Step size.}
Figures~\ref{fig:sim-bias-delta} and~\ref{fig:sim-rmse-delta} and Table~\ref{tab:sim-step} show the bias and the RMSE of $\widehat a_T$ and $\widehat b_T$ at fixed horizon $T=100$ as functions of $\delta$ (from smallest to largest step). For the largest step size $\delta=0.1$, both estimators are biased downwards, as shown in Figure~\ref{fig:sim-bias-delta}. As step $\delta$ decreases (from right to left on the plot), the average bias increases monotonically and changes sign. The smaller the H parameter, the smaller the step size at which the bias sign changes. 
One can observe that the corresponding bias does not tend to 0 as $\delta$ decreases. We investigated this behaviour, and simulations indicate that bias can be separated into two components. One is approximately linear in $\delta$ (see \eqref{eq:bias-decomposition} below), while the second one represents the bias of the continuous-time estimator at the horizon, which grows with $H$.

For $H\le0.7$ the RMSE (Figure~\ref{fig:sim-rmse-delta}) decreases as $\delta$ decreases and changes little for $\delta\le0.02$. For $H\ge0.8$ it is smallest at intermediate step sizes, $\delta\approx0.025$--$0.05$. For $H=0.9$ the RMSE of $\widehat b_T$ at $\delta=0.005$ ($0.24$) even exceeds that at $\delta=0.1$ ($0.20$).

\begin{table}%[ht]
\centering
\footnotesize
\caption{Step-size study: bias (standard error) and RMSE of $\hat a_T$ and $\hat b_T$ for $T=100$.}
\label{tab:sim-step}
\begin{tabular}{rrrrrr}
\toprule
$H$ & $\delta$ & Bias $\hat a$ (SE) & Bias $\hat b$ (SE) & RMSE $\hat a$ & RMSE $\hat b$ \\
\midrule
0.6 & 0.005 & $0.006\;(0.003)$  & $0.010\;(0.004)$  & $0.080$ & $0.120$ \\
    & 0.01  & $-0.005\;(0.003)$ & $-0.006\;(0.004)$ & $0.080$ & $0.119$ \\
    & 0.02  & $-0.026\;(0.002)$ & $-0.037\;(0.004)$ & $0.083$ & $0.122$ \\
    & 0.025 & $-0.036\;(0.002)$ & $-0.051\;(0.004)$ & $0.086$ & $0.126$ \\
    & 0.04  & $-0.064\;(0.002)$ & $-0.093\;(0.004)$ & $0.100$ & $0.146$ \\
    & 0.05  & $-0.081\;(0.002)$ & $-0.119\;(0.003)$ & $0.111$ & $0.162$ \\
    & 0.1   & $-0.154\;(0.002)$ & $-0.228\;(0.003)$ & $0.169$ & $0.249$ \\
\midrule
0.9 & 0.005 & $0.079\;(0.004)$  & $0.141\;(0.006)$  & $0.152$ & $0.240$ \\
    & 0.01  & $0.065\;(0.004)$  & $0.121\;(0.006)$  & $0.145$ & $0.226$ \\
    & 0.02  & $0.038\;(0.004)$  & $0.085\;(0.006)$  & $0.135$ & $0.202$ \\
    & 0.025 & $0.026\;(0.004)$  & $0.068\;(0.006)$  & $0.132$ & $0.192$ \\
    & 0.04  & $-0.010\;(0.004)$ & $0.017\;(0.005)$  & $0.130$ & $0.171$ \\
    & 0.05  & $-0.031\;(0.004)$ & $-0.016\;(0.005)$ & $0.132$ & $0.164$ \\
    & 0.1   & $-0.117\;(0.004)$ & $-0.149\;(0.004)$ & $0.169$ & $0.203$ \\
\bottomrule
\end{tabular}
\end{table}

\begin{figure}%[ht]
\centering
\includegraphics[width=\textwidth]{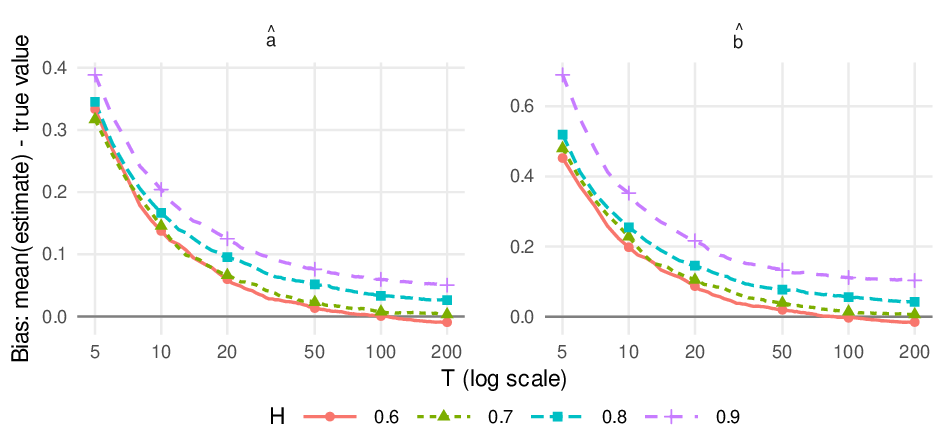}
\caption{Horizon study, $\delta=0.01$: bias of $\widehat a_T$ (left) and $\widehat b_T$ (right) as a function of $T$ (logarithmic scale), by $H$}
\label{fig:sim-bias-T}
\end{figure}

\begin{figure}%[ht]
\centering
\includegraphics[width=\textwidth]{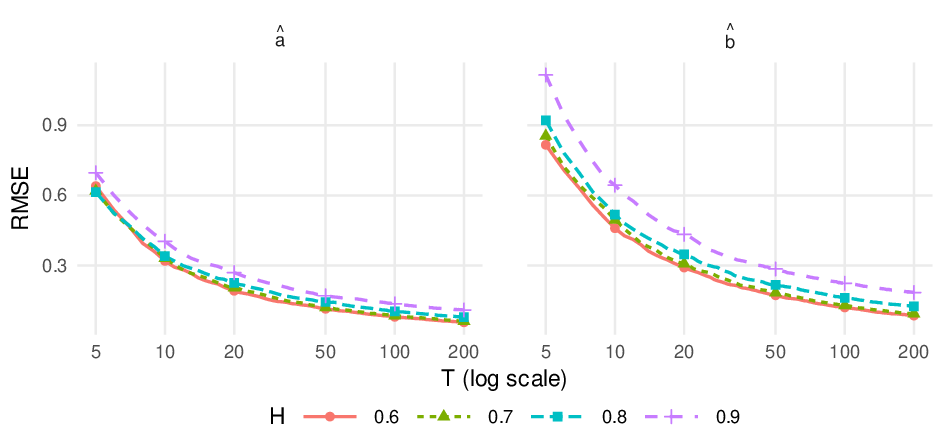}
\caption{Horizon study, $\delta=0.01$: RMSE of $\widehat a_T$ (left) and $\widehat b_T$ (right) as a function of $T$ (logarithmic scale), by $H$}
\label{fig:sim-rmse-T}
\end{figure}

\textbf{Horizon.}
Next, we used a fixed step $\delta = 0.01$ and studied the effect of increasing the horizon $T$, as shown in Figures~\ref{fig:sim-bias-T}, \ref {fig:sim-rmse-T}, and Table~\ref{tab:sim-horizon}. For smaller $T$, the bias is positive and tends to $0$ as $T$ increases. At the same time, it has a non-zero trend, which is especially clear for larger parameters $H$. Our simulations showed that, with fixed $\delta$, increasing $T$ improves accuracy up to a point, which seems to be limited by the discretization error proportional to the step size.  
At $T=200$, the bias is specifically proportional to the order of the discretization error for $H\in\{0.6,0.7\}$, whereas for $H=0.9$ it decreases slowly and is still about $0.05$ for $\widehat a_T$ and $0.1$ for $\widehat b_T$. Our simulation results also indicate that, for $H=0.6$, the RMSE decreases approximately as $T^{-1/2}$. For larger $H$ it decreases more slowly.

\begin{table}%[ht]
    \centering
    \caption{Horizon study: bias (standard error) and RMSE of $\hat a_T$ and $\hat b_T$ for $\delta=0.01$.}
    \label{tab:sim-horizon}
    \footnotesize
    \begin{tabular}{rrrrrr}
    \toprule
        $H$ & $T$ & Bias $\hat a$ (SE) & Bias $\hat b$ (SE) & RMSE $\hat a$ & RMSE $\hat b$ \\
    \midrule
        0.6 & 20  & $0.060\;(0.006)$  & $0.087\;(0.009)$  & $0.192$ & $0.291$ \\
            & 50  & $0.014\;(0.004)$  & $0.019\;(0.005)$  & $0.115$ & $0.173$ \\
            & 100 & $0.001\;(0.003)$  & $-0.003\;(0.004)$ & $0.081$ & $0.121$ \\
            & 200 & $-0.009\;(0.002)$ & $-0.015\;(0.003)$ & $0.057$ & $0.086$ \\
    \midrule
        0.7 & 20  & $0.066\;(0.006)$  & $0.105\;(0.009)$  & $0.205$ & $0.308$ \\
            & 50  & $0.023\;(0.004)$  & $0.039\;(0.006)$  & $0.121$ & $0.184$ \\
            & 100 & $0.008\;(0.003)$  & $0.014\;(0.004)$  & $0.086$ & $0.130$ \\
            & 200 & $0.003\;(0.002)$  & $0.006\;(0.003)$  & $0.062$ & $0.094$ \\
    \midrule
        0.8 & 20  & $0.095\;(0.006)$  & $0.145\;(0.010)$  & $0.226$ & $0.348$ \\
            & 50  & $0.052\;(0.004)$  & $0.077\;(0.006)$  & $0.144$ & $0.217$ \\
            & 100 & $0.033\;(0.003)$  & $0.056\;(0.005)$  & $0.104$ & $0.162$ \\
            & 200 & $0.026\;(0.002)$  & $0.042\;(0.004)$  & $0.080$ & $0.126$ \\
    \midrule
        0.9 & 20  & $0.125\;(0.008)$  & $0.216\;(0.012)$  & $0.270$ & $0.434$ \\
            & 50  & $0.076\;(0.005)$  & $0.133\;(0.008)$  & $0.171$ & $0.285$ \\
            & 100 & $0.060\;(0.004)$  & $0.112\;(0.006)$  & $0.136$ & $0.224$ \\
            & 200 & $0.050\;(0.003)$  & $0.103\;(0.005)$  & $0.110$ & $0.185$ \\
    \bottomrule
    \end{tabular}
\end{table}

\begin{figure}%[ht]
\centering
\includegraphics[width=\textwidth]{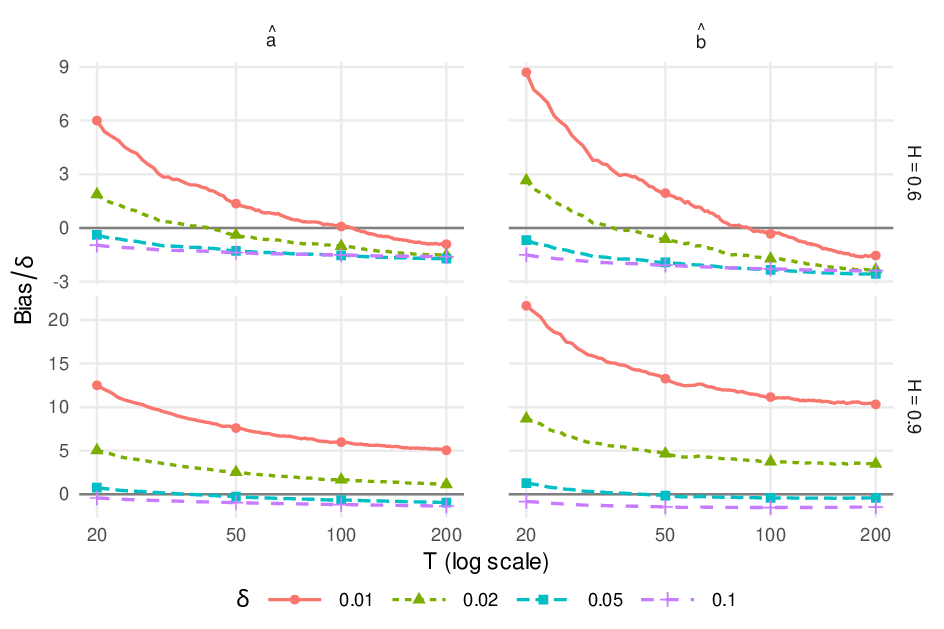}
\caption{Horizon study: bias divided by the step size, $(\widehat a_T-a)/\delta$ (left) and $(\widehat b_T-b)/\delta$ (right), as functions of $T\in[20,200]$, for $H=0.6$ (top) and $H=0.9$ (bottom); one curve per $\delta$}
\label{fig:sim-bias-over-delta}
\end{figure}

\textbf{Structure of the bias.} Investigating the results of our two studies, we came to the conclusion that our bias depends on two components: how long our data are observed (horizon $T$) and how frequently (step size $\delta$), with bias depending on both, and linearly on the second one:
\begin{equation}\label{eq:bias-decomposition}
    \widehat a_T-a\approx c_a(T)+k_a\delta,
    \qquad
    \widehat b_T-b\approx c_b(T)+k_b\delta,
\end{equation}
where $c_a(T)$, $c_b(T)$ are the horizon-specific components of the bias at the horizon $T$ and $k_a\delta$, $k_b\delta$ are the leading terms of the discretization error, with $k_a,k_b<0$. 
Dividing by $\delta$, we obtain $(\widehat a_T-a)/\delta\approx c_a(T)/\delta+k_a$, and similarly for $\widehat b_T$. Thus, if the finite-horizon bias is negligible, the curves for different $\delta$ coincide at the level $k_a$. Otherwise, they are separated in proportion to $1/\delta$.
That is exactly what we observe in Figure~\ref{fig:sim-bias-over-delta} and Table~\ref{tab:sim-bias-over-delta}. For $H=0.6$, the curves approach one specific limit as $T$ grows.
In our case, we obtained $k_a\approx-1.6$ and $k_b\approx-2.4$. For $H=0.9$, the finite-horizon bias is positive and decreases slowly, so the curves remain separated up to $T=200$.

\begin{table}%[ht]
    \centering
    \caption{Horizon study: bias divided by the step size, $(\hat a_T-a)/\delta$ and $(\hat b_T-b)/\delta$, with standard errors (divided by $\delta$) in parentheses.}
    \label{tab:sim-bias-over-delta}
    \footnotesize
    \begin{tabular}{lrrrrrr}
    \toprule
    & & & \multicolumn{4}{c}{$\delta$} \\
    \cmidrule(lr){4-7}
    Estimator & $H$ & $T$ & $0.01$ & $0.02$ & $0.05$ & $0.1$ \\
    \midrule
    $\hat a$ & 0.6 & 20  & $6.00\;(0.58)$  & $1.87\;(0.28)$  & $-0.39\;(0.11)$ & $-0.96\;(0.05)$ \\
             &     & 50  & $1.36\;(0.36)$  & $-0.40\;(0.18)$ & $-1.28\;(0.07)$ & $-1.38\;(0.03)$ \\
             &     & 100 & $0.08\;(0.25)$  & $-1.02\;(0.13)$ & $-1.53\;(0.05)$ & $-1.51\;(0.02)$ \\
             &     & 200 & $-0.91\;(0.18)$ & $-1.52\;(0.09)$ & $-1.72\;(0.03)$ & $-1.60\;(0.02)$ \\
             & 0.9 & 20  & $12.51\;(0.76)$ & $5.04\;(0.37)$  & $0.76\;(0.14)$  & $-0.43\;(0.07)$ \\
             &     & 50  & $7.60\;(0.48)$  & $2.53\;(0.24)$  & $-0.28\;(0.09)$ & $-0.97\;(0.05)$ \\
             &     & 100 & $5.98\;(0.39)$  & $1.67\;(0.19)$  & $-0.67\;(0.08)$ & $-1.18\;(0.04)$ \\
             &     & 200 & $5.04\;(0.31)$  & $1.12\;(0.16)$  & $-0.96\;(0.06)$ & $-1.36\;(0.03)$ \\
    \midrule
    $\hat b$ & 0.6 & 20  & $8.70\;(0.88)$  & $2.66\;(0.43)$  & $-0.68\;(0.16)$ & $-1.50\;(0.07)$ \\
             &     & 50  & $1.94\;(0.54)$  & $-0.63\;(0.27)$ & $-1.93\;(0.10)$ & $-2.09\;(0.05)$ \\
             &     & 100 & $-0.33\;(0.38)$ & $-1.73\;(0.19)$ & $-2.35\;(0.07)$ & $-2.29\;(0.03)$ \\
             &     & 200 & $-1.54\;(0.27)$ & $-2.33\;(0.13)$ & $-2.57\;(0.05)$ & $-2.39\;(0.02)$ \\
             & 0.9 & 20  & $21.64\;(1.19)$ & $8.68\;(0.57)$  & $1.28\;(0.21)$  & $-0.84\;(0.09)$ \\
             &     & 50  & $13.27\;(0.80)$ & $4.65\;(0.38)$  & $-0.16\;(0.14)$ & $-1.45\;(0.06)$ \\
             &     & 100 & $11.16\;(0.61)$ & $3.73\;(0.29)$  & $-0.42\;(0.11)$ & $-1.52\;(0.05)$ \\
             &     & 200 & $10.32\;(0.48)$ & $3.47\;(0.23)$  & $-0.43\;(0.08)$ & $-1.47\;(0.04)$ \\
    \bottomrule
    \end{tabular}
\end{table}

\subsection{Ratio estimator: simulation design}
\label{sec:sim-ratio-design}
Theorem~\ref{thm:parameter-relation} applies to a broad class of Gaussian processes. In our simulations, we apply it for the two:
$G=\sigma B^H$ and $G=\sigma W^H$, that is, for model \eqref{eq:model} driven by Riemann–Liouville fractional Brownian motion and by fractional Brownian motion (see Subsection~\ref{subsec:ratio-examples}). We take the same initial parameters $\sigma=0.5$, $X_0=1.5$, but only $H\in\{0.6,0.8\}$. Additionally, we consider $b=1$ as known, so that the ratio estimator yields the estimator
\begin{equation*}
    \widehat a_T=b\,\widehat\theta_T
\end{equation*}
of $a$. In our simulation of this specific model, we consider cases of $a\in\{4,$ $2,$ $1,$ $0.5,$ $0.25\}$. The integrals in $\widehat\theta_T$ are approximated by left Riemann sums:
\begin{equation*}
\widehat\theta_T \approx \frac{\sum_{k=0}^{N-1}X_{t_k}}{\sum_{k=0}^{N-1}X_{t_k}^{-1}},
\end{equation*}
with the common factor $\delta$ canceling. Process $X$ here is constructed by the scheme of Subsection~\ref{sec:sim-ab-design}. Both processes $B^H$ and $W^H$ are simulated
exactly at the grid points using the Cholesky method described in Subsection~\ref{sec:sim-ab-design}: for $B^H$ with the covariance from Definition~\ref{def:rlfbm}, and for $W^H$ with $\cov(W^H_j,W^H_k)=\frac{1}{2}\bigl(j^{2H}+k^{2H}-|j-k|^{2H}\bigr)$ on the integer grid, scaled by $\delta^{2H}$.

Since $\widehat\theta_T$ involves no numerical approximation of fractional derivatives, its discretization error is smaller than that of $(\widehat a_T,\widehat b_T)$ (see Subsection~\ref{sec:sim-ratio-results}). 
That is why we used larger step sizes $\delta$ and proportionally longer horizons $T$, while the finest grid has $N=10^4$.
We conduct two studies, each with $m=1000$ replications. In the \emph{horizon study}, $\delta=0.1$ and the estimates are recorded at $T\in\{5,$ $6,$ $8,$ $10,$ $12,$ $15,$ $20,$ $25,$ $30,$ $40,$ $50,$ $60,$ $80,$ $100,$ $120,$ $150,$ $200,$ $250,$ $300,$ $400,$ $500,$ $600,$ $800,$ $1000\}$. In the \emph{step-size study}, $T=100$ and $\delta\in\{0.01,$ $0.02,$ $0.05,$ $0.1\}$. As for $(\widehat a_T,\widehat b_T)$, all step sizes within one replication are obtained from one trajectory of the noise on the finest grid, and all horizons from one trajectory on the longest interval.

\subsection{Ratio estimator: results}
\label{sec:sim-ratio-results}
Since the true values of $a$ differ by a factor of $16$, we report the accuracy measures of $\widehat a_T$ relative to $a$, in percent.
For $m$ replications $\widehat a_T^{(1)},\dots,\widehat a_T^{(m)}$ at a fixed design point, let $\overline a_T=\frac1m\sum_{i=1}^m\widehat a_T^{(i)}$ and $s_a^2=\frac1{m-1}\sum_{i=1}^m\bigl(\widehat a_T^{(i)}-\overline a_T\bigr)^2$. 
We define
\begin{gather*}
\mathrm{RelBias}(\widehat a_T)=\frac{\overline a_T-a}{a}\cdot100\%,\qquad
\mathrm{RelSE}(\widehat a_T)=\frac{s_a}{a\sqrt m}\cdot100\%,
\\
\mathrm{RelRMSE}(\widehat a_T)=\frac1a\biggl(\frac1m\sum_{i=1}^m
\bigl(\widehat a_T^{(i)}-a\bigr)^2\biggr)^{1/2}\cdot100\%.
\end{gather*}

\textbf{Step size.}
We illustrate the relative bias and the relative RMSE of $\widehat a_T$ as functions of $\delta$ for $T=100$ in Table~\ref{tab:sim-ratio-step}. 
Our simulations indicated that the step size $\delta$ matters less for this estimator than for $(\widehat a_T,\widehat b_T)$. Accordingly, increasing $\delta$ from $0.01$ to $0.1$ changes the relative bias by at most $0.3$ percentage points and the relative RMSE by at most $0.2$ percentage points. Changes in the relative bias are smaller than its standard errors.

\begin{figure}%[ht]
\centering
\includegraphics[width=\textwidth]{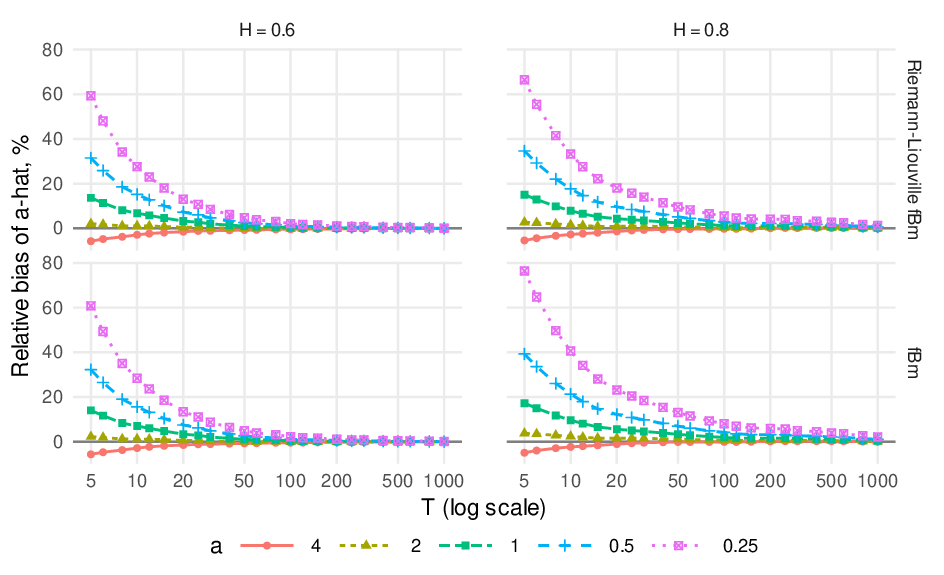}
\caption{Horizon study, $\delta=0.1$: relative bias $(\widehat a_T-a)/a$ of $\widehat a_T$, in percent, as a function of $T$ (logarithmic scale), by $a$. Rows: Riemann--Liouville fBm and fBm; columns: $H$}
\label{fig:sim-ratio-bias}
\end{figure}

\begin{figure}%[ht]
\centering
\includegraphics[width=\textwidth]{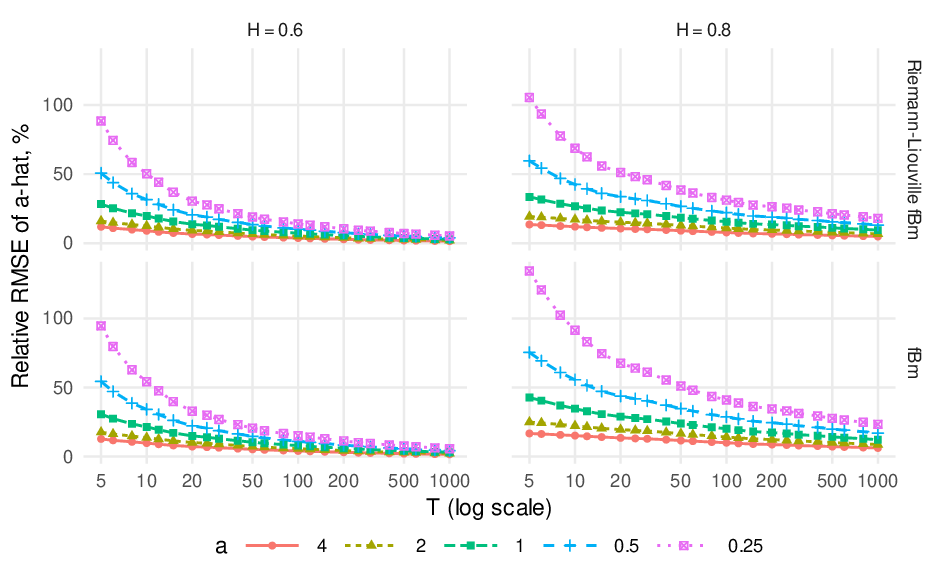}
\caption{Horizon study, $\delta=0.1$: relative RMSE of $\widehat a_T$, in percent, as a function of $T$ (logarithmic scale); panels as in Figure~\ref{fig:sim-ratio-bias}}
\label{fig:sim-ratio-rmse}
\end{figure}

\textbf{Horizon}. Figures~\ref{fig:sim-ratio-bias} and~\ref{fig:sim-ratio-rmse} show the relative bias and the relative RMSE of $\widehat a_T$ for $\delta=0.1$ as functions of $T$, while Table~\ref{tab:sim-ratio} summarizes them at $T=100$ and $T=1000$. For $T\le20$ the bias is positive for $a\le2$ and negative for $a=4$, and it is larger for smaller $a$. 
For $H=0.6$ the bias is almost the same for both noises, whereas for $H=0.8$ it is larger for $W^H$. The relative RMSE decreases monotonically in $T$ in all configurations. It is larger for smaller $a$, for $H=0.8$ than for $H=0.6$, and for $W^H$ than for $B^H$. 

\begin{table}%[ht]
\centering
\caption{Ratio estimator $\hat a_T=b\,\widehat\theta_T$ with known $b=1$, step-size study, $T=100$, $a\in\{4,1,0.25\}$: relative bias $(\hat a_T-a)/a$ (standard error) and relative RMSE, in percent.}
\label{tab:sim-ratio-step}
\footnotesize
\begin{tabular}{rrrrrrr}
\toprule
& & & \multicolumn{2}{c}{Riemann--Liouville fBm} & \multicolumn{2}{c}{fBm} \\
\cmidrule(lr){4-5}\cmidrule(lr){6-7}
&&& Rel.\ bias, & RMSE & Rel.\ bias, & RMSE\\
$H$ & $a$ & $\delta$ &\% (SE) & Rel.\, \% & \% (SE) & Rel.\, \% \\
\midrule
0.6 & 4    & 0.01 & $-0.17\;(0.11)$ & $3.57$  & $-0.16\;(0.12)$ & $3.90$ \\
    &      & 0.1  & $-0.23\;(0.11)$ & $3.57$  & $-0.21\;(0.12)$ & $3.90$ \\
    & 1    & 0.01 & $0.77\;(0.22)$  & $7.10$  & $0.82\;(0.24)$  & $7.75$ \\
    &      & 0.1  & $0.85\;(0.22)$  & $7.12$  & $0.90\;(0.24)$  & $7.78$ \\
    & 0.25 & 0.01 & $2.78\;(0.43)$  & $13.73$ & $2.93\;(0.46)$  & $14.85$ \\
    &      & 0.1  & $3.03\;(0.43)$  & $13.90$ & $3.18\;(0.47)$  & $15.04$ \\
\midrule
0.8 & 4    & 0.01 & $0.13\;(0.24)$ & $7.62$  & $0.29\;(0.31)$ & $9.68$ \\
    &      & 0.1  & $0.08\;(0.24)$ & $7.61$  & $0.24\;(0.31)$ & $9.67$ \\
    & 1    & 0.01 & $1.82\;(0.48)$ & $15.37$ & $2.47\;(0.62)$ & $19.65$ \\
    &      & 0.1  & $1.90\;(0.48)$ & $15.39$ & $2.54\;(0.62)$ & $19.68$ \\
    & 0.25 & 0.01 & $6.40\;(0.96)$ & $30.89$ & $8.82\;(1.24)$ & $40.31$ \\
    &      & 0.1  & $6.65\;(0.96)$ & $31.07$ & $9.08\;(1.25)$ & $40.54$ \\
\bottomrule
\end{tabular}
\end{table}

\begin{table}%[ht]
\centering
\footnotesize
\caption{Ratio estimator $\hat a_T=b\,\widehat\theta_T$ with known $b=1$, horizon study, $\delta=0.1$, $a\in\{4,1,0.25\}$: relative bias $(\hat a_T-a)/a$ (standard error) and relative RMSE, in percent.}
\label{tab:sim-ratio}
\begin{tabular}{rrrrrrr}
\toprule
    & & & \multicolumn{2}{c}{Riemann--Liouville fBm} & \multicolumn{2}{c}{fBm} \\
    \cmidrule(lr){4-5}\cmidrule(lr){6-7}
    & & & Rel.\ bias, &  RMSE  & Rel.\ bias, & RMSE \\
    $H$ & $a$ & $T$ & \% (SE) & Rel., \% & \% (SE) & Rel., \% \\
    \midrule
    0.6 & 4    & 100  & $-0.48\;(0.11)$ & $3.67$  & $-0.49\;(0.13)$ & $4.00$  \\
        &      & 1000 & $-0.09\;(0.04)$ & $1.42$  & $-0.10\;(0.05)$ & $1.56$  \\
        & 1    & 100  & $0.36\;(0.23)$  & $7.21$  & $0.36\;(0.25)$  & $7.87$  \\
        &      & 1000 & $-0.05\;(0.09)$ & $2.79$  & $-0.06\;(0.10)$ & $3.04$  \\
        & 0.25 & 100  & $2.12\;(0.44)$  & $13.93$ & $2.17\;(0.47)$  & $15.07$ \\
        &      & 1000 & $0.07\;(0.17)$  & $5.23$  & $0.07\;(0.18)$  & $5.66$  \\
    \midrule
    0.8 & 4    & 100  & $-0.23\;(0.24)$ & $7.74$  & $-0.10\;(0.31)$ & $9.95$  \\
        &      & 1000 & $-0.06\;(0.15)$ & $4.78$  & $-0.02\;(0.19)$ & $6.14$  \\
        & 1    & 100  & $1.30\;(0.49)$  & $15.55$ & $1.92\;(0.63)$  & $20.10$ \\
        &      & 1000 & $0.22\;(0.30)$  & $9.42$  & $0.43\;(0.38)$  & $12.13$ \\
        & 0.25 & 100  & $5.54\;(0.97)$  & $31.15$ & $8.03\;(1.27)$  & $40.97$ \\
        &      & 1000 & $1.26\;(0.57)$  & $17.93$ & $2.13\;(0.73)$  & $23.29$ \\
\bottomrule
\end{tabular}
\end{table}

\FloatBarrier

\end{document}